\documentclass{article}
\usepackage[nonatbib,final]{neurips_2021}
\usepackage[utf8]{inputenc}
\usepackage[T1]{fontenc}   
\usepackage{hyperref}      
\usepackage{url}          
\usepackage{booktabs}     
\usepackage{amsfonts}    
\usepackage{nicefrac}      
\usepackage{microtype} 
\usepackage{xcolor}    
\usepackage{lipsum}
\usepackage{acro}
\usepackage[]{algorithm2e}
\usepackage{amsmath}
\usepackage{amssymb}
\usepackage{amsthm}
\usepackage{cleveref}
\usepackage{color}
\usepackage{comment}
\usepackage{mathtools}
\usepackage[sort&compress,numbers]{natbib}
\usepackage{url}
\usepackage{cancel} 
\usepackage{graphicx}
\usepackage{lipsum}
\usepackage{float}
\usepackage{subfig}
\usepackage{placeins}
\usepackage[normalem]{ulem}
\usepackage{todonotes}
\usepackage{bm} 
\usepackage{makecell} 
\usepackage{nicematrix}  

\allowdisplaybreaks

\acsetup{first-long-format=\itshape}
\DeclareAcronym{GP}{short=GP, long=Gaussian process, short-plural=s, long-plural=es}
\DeclareAcronym{BBPN}{short=BBPN, long=black box probabilistic numerics}
\DeclareAcronym{RDAL}{short=RDAL, long=Richardson's deferred approach to the limit}
\DeclareAcronym{PN}{short=PN, long=probabilistic numerics}
\DeclareAcronym{PDE}{short=PDE, long=partial differential equation}
\DeclareAcronym{ODE}{short=ODE, long=ordinary differential equation}
\DeclareAcronym{IVP}{short=IVP, long=initial value problem}
\DeclareAcronym{RKHS}{short=RKHS, long=reproducing kernel Hilbert space}
\DeclareAcronym{KSE}{short=KSE, long=Kuramoto--Sivashinsky equation}

\newtheorem{definition}{Definition}

\newtheorem{remark}[definition]{Remark}
\newtheorem*{remark*}{Remark}

\newtheorem{proposition}[definition]{Proposition}
\newtheorem{example[definition]}{Example}
\newtheorem*{example*}{Example}

\newcommand{\indep}{\perp \!\!\! \perp}
\def\d{\,\mathrm{d}}

\usepackage{everypage}
\usepackage{afterpage}
\usepackage{ifthen}
\usepackage[width=0.89\textwidth]{caption}
\usepackage{refcount} 

\title{Black Box Probabilistic Numerics}

\author{%
  Onur Teymur \\
  University of Kent\\
  Alan Turing Institute \\
  \And
  Christopher N. Foley \\
  University of Cambridge \\
  Optima Partners
  \And
  Philip G. Breen \\
  Roar AI \\
  \AND
  Toni Karvonen \\
  University of Helsinki \\
  Alan Turing Institute
  \And
  Chris. J. Oates \\
  Newcastle University \\
  Alan Turing Institute \\
}

\begin{document}

\maketitle

\begin{abstract}
\emph{Probabilistic numerics} casts numerical tasks, such the numerical solution of differential equations, as inference problems to be solved.
One approach is to model the unknown quantity of interest as a random variable, and to constrain this variable using data generated during the course of a traditional numerical method.
However, data may be nonlinearly related to the quantity of interest, rendering the proper conditioning of random variables difficult and limiting the range of numerical tasks that can be addressed. 
Instead, this paper proposes to construct probabilistic numerical methods based only on the final output from a traditional method.
A convergent sequence of approximations to the quantity of interest constitute a dataset, from which the limiting quantity of interest can be extrapolated, in a probabilistic analogue of Richardson's deferred approach to the limit.
This \emph{black box} approach (1) massively expands the range of tasks to which probabilistic numerics can be applied, (2) inherits the features and performance of state-of-the-art numerical methods, and (3) enables provably higher orders of convergence to be achieved.
Applications are presented for nonlinear ordinary and partial differential equations, as well as for eigenvalue problems\textemdash a setting for which no probabilistic numerical methods have yet been developed.
\end{abstract}

\section{Introduction} \label{sec: intro}
\Ac{PN} has attracted significant recent interest from researchers in machine learning, motivated by the possibility of incorporating probabilistic descriptions of \emph{numerical uncertainty} into applications of probabilistic inference and decision support \cite{hennig15,oates2019modern}. 
\ac{PN} treats the intermediate calculations performed in running a traditional (\emph{i.e.} non-probabilistic)
numerical procedure as \emph{data}, which can be used to constrain a random variable model for the quantity of interest \citep{cockayne2019bayesian}. Conjugate Gaussian inference has been widely exploited, with an arsenal of \ac{PN} methods developed for linear algebra \citep{cockayne2019abayesian,bartels2019probabilistic,wenger2020probabilistic,hennig2015probabilistic,reid2020probabilistic,schafer2017compression,bartels2016probabilistic,cockayne2020probabilistic}, cubature \citep{diaconis88,ohagan92,fisher2020locally,pruher2016use,gessner2020active,karvonen2019symmetry,chai2019improving,jagadeeswaran2019fast,karvonen2017classical,karvonen2018bayes,osborne2012bayesian,xi2018bayesian,briol2015frank,gunter2014sampling,kennedy1998bayesian,o1991bayes,larkin1972gaussian,rasmussen2003bayesian,briol2019probabilistic}, optimisation \citep{mockus77,mockus78,mockus89,snoek12,hennig13a,mahsereci15}, and differential equations \citep{skilling92,conrad16,chkrebtii16,schober18,teymur16,teymur18,hennig13,kersting16,schober14,owhadi2015bayesian,tronarp18,wang18,cockayne17,chkrebtii19,owhadi2019operator,abdulle20,kersting20,wang2021bayesian,bosch21}.
However, nonlinear tasks pose a major technical challenge to this approach, as well as to computational statistics in general, due to the absence of explicit conditioning formulae.
Compared to traditional numerical methods, which have benefited from a century or more of sustained research effort, the current scope of \ac{PN} is limited. 
The performance gap is broadly characterised by the absence of certain important functionalities---adaptivity, numerical well-conditioning, efficient use of computational resource---all of which contribute to limited applicability in real-world settings.

This article proposes a pragmatic solution that enables state-of-the-art numerical algorithms to be immediately exploited in the context of \ac{PN}.
The idea, which we term \ac{BBPN}, is a statistical perspective on \ac{RDAL} \cite{richardson27}.
The starting point for \ac{BBPN} is a sequence of increasingly accurate approximations produced by a traditional numerical method as its computational budget is increased. 
Extrapolation of this sequence (to the unattainable limit of `infinite computational budget') is formulated as a prediction task, to which statistical techniques can be applied. 
For concreteness, we perform this prediction using \acp{GP} \cite{rasmussen06}, but other models could be used.
Note that we do not aim to remove the numerical analyst from the loop; the performance of \ac{BBPN} is limited by that of the numerical method on which it is based.

There are three main advantages of \ac{BBPN} compared to existing methods in \ac{PN}: (1) \ac{BBPN} is applicable to any numerical task for which there exists a traditional numerical method; (2) state-of-the-art performance and functionality are automatically inherited from the underlying numerical method; (3) \ac{BBPN} achieves a provably higher order of convergence relative to a single application of the numerical method on which it is based, in an analogous manner to \ac{RDAL}. 
The main limitations of \ac{BBPN}, compared to existing methods in \ac{PN}, are: (1) multiple realisations of a traditional numerical method are required (\emph{i.e.} one datum is not sufficient in an extrapolation task), and (2) a joint statistical model has to be built for not just the quantity of interest (as in standard \ac{PN}), but also for the error associated with the output of a traditional numerical method.
The capacity of generic statistical models, such as \acp{GP}, to learn salient aspects of this structure from data and to produce meaningful predictions over a range of real-world numerical problems, demands to be investigated.

The article is organised as follows:
In \Cref{sec: rdal} we recall classical \ac{RDAL}.
In \Cref {sec:setting} we lift \ac{RDAL} to the space of probability distributions, exploiting \acp{GP} to instantiate \ac{BBPN} and providing a theoretical guarantee that higher order convergence is achieved by using \acp{GP} within \ac{BBPN}.
In \Cref{sec: empirical} we present a detailed empirical investigation into \ac{BBPN}, demonstrating its effectiveness on challenging tasks that go beyond the capability of current methods in \ac{PN}, while also highlighting potential pitfalls. As part of this we perform a comparison of the uncertainty quantification properties of \ac{BBPN} against earlier approaches.
A closing discussion is contained in \Cref{sec: discuss}.

\section{Turning Lead into Gold} \label{sec: rdal}

Our starting point is the celebrated observation of Richardson \citep{richardson27}, that multiple numerical approximations can be combined to produce an approximation more accurate than any of the individual approximations.
To see this, consider an intractable scalar quantity of interest $q^\ast \in \mathbb{R}$, 
and suppose that $q^*$ can be approximated by a numerical method $q$ that depends on a parameter $h>0$, such that 
\begin{equation}
q(h) = q^\ast + Ch^\alpha + \mathcal{O}(h^{\alpha+1}) \label{eq: def of convergent nm}
\end{equation}
for some $C\in \mathbb{R}$ (which may be unknown) and $\alpha > 0$ (which is assumed known, and called the \emph{order} of the method).
Clearly $q(h)$ converges to $q^\ast$ as $h \rightarrow 0$, but we also suppose that the \emph{cost} of computing $q(h)$ increases in the same limit, with exact evaluation of $q(0)$ requiring a hypothetically infinite computational budget.
\Cref{prop: richard} is the cornerstone of \ac{RDAL}. 
It demonstrates that two evaluations of a numerical method of order $\alpha$ can be combined to obtain a numerical method of order $\alpha + 1$. An elementary proof of this foundational result is provided in \Cref{ap: richard proof}.

\begin{proposition} \label{prop: richard}
Let $q$ be a numerical method of order $\alpha$, as in \eqref{eq: def of convergent nm}.
Fix $\gamma \in (0,1)$  and let $q_\gamma(h)$ denote the height at which a straight line drawn through the points $(h^\alpha,q(h))$ and $((\gamma h)^\alpha,q(\gamma h))$ intersects the vertical axis in $\mathbb{R}^2$.
Then $q_\gamma$ is a numerical method of order $\alpha+1$.
\end{proposition}

Now consider the natural generalisation of \Cref{prop: richard}, in which we compute approximations $q(h_i)$ along a decreasing sequence of values $(h_i)_{i=1}^n$.  
One can then fit a smooth interpolant to the points $\{(h_i^\alpha,q(h_i))\}_{i=1}^n$ (generalising the straight line through two points), then extrapolate this to $h = 0$, to give an estimate for the quantity of interest. 
This simple idea is widely used in numerical analysis; its potential to radically improve solution accuracy, given only a sequence of simple calculations as input, prompted Press et al. \cite[p.\@~922]{press07} to describe it as ``turning lead into gold''.
The practical success of \ac{RDAL} depends on the choice of interpolant, with polynomial interpolation being most commonly used.
Unqualified, \ac{RDAL} is usually understood to refer to an order $n-1$ polynomial fitted to $n$ points, which produces a numerical method of order $\alpha+n$; see Theorem 9.1 of \cite{hairer08}. 
Higher-order polynomial extrapolation is known to perform poorly unless the values $(h_i)_{i=1}^n$ are able to be chosen specifically to mitigate Runge's phenomenon \citep{runge1901empirische}, motivating the Bulirsch--Stoer algorithm \cite{bulirsch64}, which instead fits a rational function interpolant. This allows both greater expressiveness and robustness than polynomial interpolation (though not necessarily as efficiently \cite{press07}). These methods are all situated within the broad category of \emph{extrapolation methods} in numerical analysis; a comprehensive historical survey can be found in \cite{joyce71}.

\begin{figure}[!t]
\vspace*{-1em}
    \centering
   \hspace*{-1em} \includegraphics[width=1.03\textwidth]{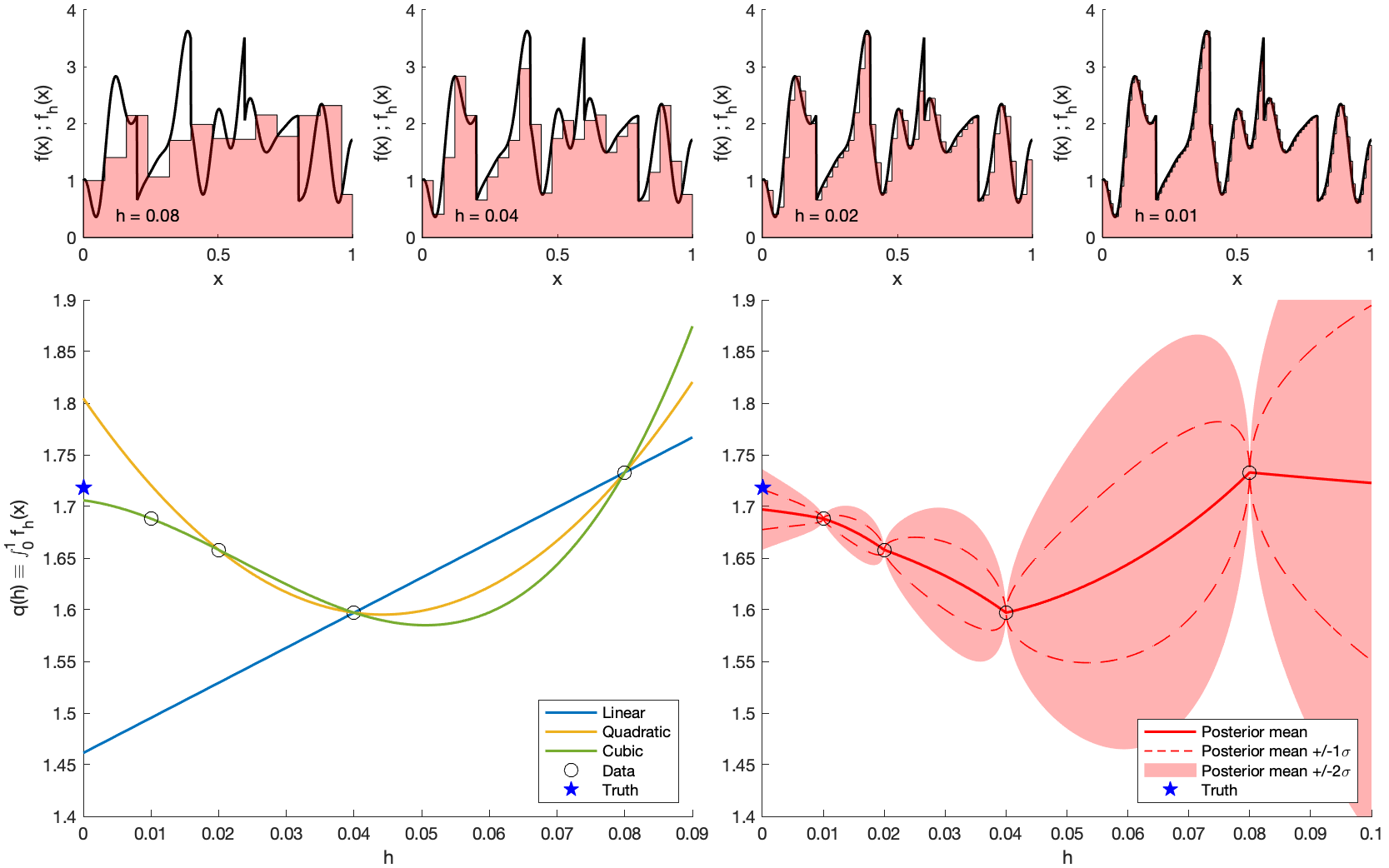}
    \caption{Richardson's deferred approach to the limit (RDAL), applied to the method of Riemann sums with integration bandwidth $h$ and an oscillatory integrand $f(x)$, displayed in the four top panes.
    The bottom left pane shows linear, quadratic and cubic interpolants converging on the true value of the integral, denoted by a blue star. 
    The strength of RDAL is seen in the fact that the cubic interpolant gives a better estimate than that given by the finest-grid Riemann sum with $h=0.01$. (This example is closely related to Romberg's method.)
    The bottom right pane illustrates \emph{black box probabilistic numerics} (BBPN), in which a Gaussian process is fitted to the same data.
    The GP specification is crucial to the performance of BBPN, and is described in \Cref{sec:setting}.
    }
    \label{fig:classical example}
    \vspace*{-0.5em}
\end{figure}

\Cref{fig:classical example} presents a simple visual demonstration of \ac{RDAL}, applied to the method of Riemann sums for an oscillatory 1D integrand.
While \ac{RDAL} gives improved approximations, no quantification of estimation uncertainty is provided.
The only attempt of which we are aware to provide uncertainty quantification for \ac{RDAL} is due to \cite{oliver2014estimating}, who focused on the Navier--Stokes equation and a specific scalar quantity of interest.
Here we go further, proposing the general framework of \ac{BBPN} and introducing novel methodology that goes beyond scalar quantities of interest.
The right-hand pane of \Cref{fig:classical example} displays the outcome of the method we are about to introduce, applied to the same task\textemdash observe that the true value of the integral falls within the $\pm 2 \sigma$ credible set produced using \ac{BBPN}.
Details of the simulations in this figure are contained in \Cref{app: Riemann}.

\section{Methodology} \label{sec:setting}
The core idea of \ac{BBPN} is to model $q(h)$ as a \emph{stochastic process} $Q(h)$ rather than fit a deterministic interpolant as in \ac{RDAL}.
The distribution of the marginal random variable $Q(0)$ is then interpreted as a representation of the epistemic uncertainty in the quantity of interest $q(0)$.
Conjugate Gaussian inference can be performed in \ac{BBPN}, since one needs only to construct an interpolant. 
This means the challenge of nonlinear conditioning encountered in \ac{PN} \citep{cockayne2019abayesian} is avoided, massively extending the applicability of \ac{PN}.
In addition to being able to leverage state-of-the-art numerical methods, the \ac{BBPN} approach enjoys provably higher orders of convergence relative to a single application of the numerical method on which it is based; see \Cref{subsec: GPs}. 
\subsection{Notation and Setup}
Our starting point is to generalise \eqref{eq: def of convergent nm} to encompass essentially all numerical tasks, following the abstract perspective of \cite{chartres1972general}.
To do so, we observe that any quantity of interest $q^*$ can be characterised by a sufficiently large collection of real values $q^*(t)$, with $t$ ranging over an appropriate index set $T$.

\begin{definition} \label{def: tnm}
A traditional numerical method is defined as a map $q : [0,h_0) \times T \rightarrow \mathbb{R}$, for some $h_0 > 0$ such that, for all $t \in T$, the function $h \mapsto q(h,t)$ is continuous at $0$ with limit $q(0,t) = q^*(t)$. 
\end{definition}

\noindent For example, a (univariate) \ac{IVP}, in which $t$ is interpreted as time, can be solved using a traditional numerical method $q$ whose time step size $h$ trades off approximation error against computational cost.
The output $q(h,t)$ of such a method represents an approximation to the true solution $q^*(t)$ of the \ac{IVP}, at each time $t$ for which the solution is defined.
In general, depending on the numerical task, the index $t$ could be spatio-temporal, discrete, or even an unstructured set, while the meaning of the index $h$ will depend on the numerical method.

\Cref{def: tnm} thus encompasses, among other things: (1) adaptive integrators for time-evolving \acp{PDE}, where $h > 0$ represents a user-specified error tolerance, and the spatio-temporal domain of the solution is indexed by $T$; (2) iterative methods for approximating the singular values of a $d \times d$ matrix, where for example $h := w^{-1}$ with $w$ the number of iterations performed, and the ordered singular values are indexed by $T = \{1,\dots,d\}$; and (3) the simultaneous approximation of multiple related quantities of interest, where $T$ indexes not only the domain(s) on which the individual quantities of interest are defined, but also the multiple quantities of interest themselves (this situation arises, for example, in inverse problems that are \textit{\ac{PDE}-constrained} \citep{de2015numerical}).

The perspective in \Cref{def: tnm} is abstract but, as these examples make clear, it will typically only be possible to compute $q$ at certain input values (such as $h = w^{-1}$ for $w \in \mathbb{N}$, or for just a finite collection of inputs $t$ if the index set $T$ is infinite), and furthermore each evaluation is likely to be associated with a computational cost.
Thus complete information regarding the map $q : [0,h_0) \times T \rightarrow \mathbb{R}$ will not be available in general, and there will therefore remain \emph{epistemic uncertainty} in its complete description.
Our aim in \Cref{subsec: bbpn sub}, in line with the central philosophical perspective of \ac{PN}, is to characterise this uncertainty using a statistical model.

\subsection{Black Box Probabilistic Numerics} \label{subsec: bbpn sub}

The proposed \ac{BBPN} approach begins with a \textit{prior} stochastic process $Q$ and constrains this prior using data $D$.
Concretely, we assume that the real values $q(h_i,t_{i,j})$ are provided at a finite set of resolutions $h_1 > \dots > h_n > 0$ and distinct ordinates $t_{i,1},\dots,t_{i,m_i} \in T$. 
Note that the number of $t$-ordinates $m_i$ can depend on $h_i$.
Our dataset therefore contains the following information on $q$:
\begin{equation}
D := \{(h_i,t_{i,j},q(h_i,t_{i,j})): \; i = 1,\dots,n; \; j = 1,\dots,m_i\} \label{eq: triples of data}
\end{equation}
The stochastic process obtained by conditioning $Q$ on the dataset $D$, denoted $Q | D$, implies a marginal distribution for $Q(0,\cdot)$, which we interpret as a statistical prediction for the unknown quantity of interest $q^\ast(\cdot)$.
In order for uncertainty quantification in this model to be meaningful, one either requires expert knowledge about the numerical method that generated $D$, or one must employ a stochastic process that is able to adapt to the data, so that its predictions can be calibrated.

Our goal is to specify a stochastic process model $Q(h,t)$ that behaves in a desirable way under extrapolation to $h = 0$.
To this end, we decompose
\begin{align}
    Q(h,t) & = Q^*(t) + E(h,t)   \label{eq: high level model}
\end{align}
where $Q^*(t)$ is a prior model for the unknown quantity of interest $q^*(t)$, and $E(h,t)$ is a prior model for the error of the numerical method.
It will be assumed that $Q^*$ and $E$ are independent (denoted $Q^* \indep E$), meaning that prior belief about the quantity of interest is independent of prior belief regarding the performance of the numerical method.
(This assumption is made only to simplify the model specification, but if detailed insight into the error structure of a numerical method is available then this can be exploited.)
Compared to the existing \ac{PN} methods cited in \Cref{sec: intro}, a prior model for the error $E$ is an additional requirement in \ac{BBPN}.

The error $E(h,t)$ is assumed to vanish\footnote{This statement covers several potentially subtle notions from numerical analysis such as well-posedness of the problem and numerical stability of the algorithm; these are studied in detail in their own right in the literature, and for our purposes it suffices to assume that the error behaves well in the limit.} as $h \rightarrow 0$, meaning that a stationary stochastic process model for $E(h,t)$, and hence for $Q(h,t)$, is inappropriate, and can result in predictions that are both severely biased as well as under-confident; see \Cref{app: prior sensitivity}.
In the next section, we propose a parsimonious non-stationary \ac{GP} model for $Q(h,t)$ of the form \eqref{eq: high level model}, which combines knowledge of the \textit{order} of the numerical method (only) with data-driven estimation of \ac{GP} hyperparameters.
This setting is practically relevant---the order of a numerical method is typically one of the first theoretical properties that researchers aim to establish while, conversely, for more complex numerical methods the order may actually be the only salient high-level error-characterising property that is known, and thus represent the limit of mathematical insight into the method.

\subsection{Gaussian Process \ac{BBPN}}
\label{subsec: GPs}

Gaussian processes provide a convenient model for $Q^*$ and $E$, since they produce an explicit form for the conditional $Q|D$.
The details of conjugate Gaussian inference are standard (see e.g. \citep{rasmussen06}) and so relegated to \Cref{app: conditioning formulae}; our focus here is on the specification of \ac{GP} priors for $Q^*$ and $E$.

The notation $Q \sim \mathcal{GP}(\mu_Q,k_Q)$ will be used to denote that $Q$ is a \ac{GP} with mean function $\mu_Q$ and covariance function $k_Q$.
With no loss of generality, in what follows we consider centred processes (\emph{i.e.} $\mu_Q = 0$).
It will be assumed that $T = T_1 \times \dots \times T_p$, with each $T_i$ either a discrete or a continuous subset of a Euclidean space, with the Euclidean distance between elements $t_i,t_i' \in T_i$ being denoted $\|t_i - t_i'\|$.
(Typical applications involve small $p$; for example, the domain of a spatio-temporal \ac{PDE} is typically decomposed as $T = T_1 \times T_2$ where $T_1$ indexes time and $T_2$ indexes all spatial dimensions, so that $p=2$.) 

\paragraph{Prior for $Q^*$:}
In the absence of detailed prior belief about $q^*$, we consider the following default prior model.
Let $G \sim \mathcal{GP}(0,\sigma^2 \rho_G k_G) $, $Z = (Z_1, \dots, Z_v) \sim \mathcal{N}(0,\sigma^2 I)$, and let $Z \indep G$. 
Let $b_1,\dots,b_v$ be a finite a collection of basis functions and set $b(t) = (b_1(t),\dots,b_v(t))^\top$.
Then set
\begin{align*}
Q^*(t) & = Z \cdot b(t) + G(t)
\end{align*}
where $\sigma^2, \rho_G > 0$ are parameters to be estimated.
The basis $b$ will be problem-specific and could be a polynomial basis, Fourier basis, or any number of other bases depending on context.
The case $v=1$ with a constant intercept is closely related to \emph{ordinary kriging} and the case $v > 1$ is closely related to \emph{universal kriging} \citep[p.\@~8]{stein2012interpolation}.
The apparent redundancy in the parameterisation due to the product $\sigma^2 \rho_G$ will be explained later.
Using the notation $t = (t_1,\dots,t_p)$, we consider a tensor product covariance model
$k_G(t,t') = \prod_{i=1}^p k_{G,i}(t_i,t_i')$, $k_{G,i}(t_i,t_i') = \phi_i\left( \|t_i-t_i'\| / \ell_{t,i} \right)$,
for some radial basis functions $\phi_i$, scaled to satisfy $\phi_i(0) = 1$, and length-scale parameters $\ell_{t,i} > 0$ to be estimated.

\paragraph{Prior for $E$:} 
The process $E(h,t)$ is a model for the numerical error $q(h,t) - q^*(t)$, $t>0$, which may be highly structured. A flexible prior model is therefore required.
Moreover the error will, by definition, depend on the order of the numerical method; for successful extrapolation we must therefore encode this order into the model for $E$.
It was earlier argued that a stationary \ac{GP} is inappropriate, since the error is assumed to be $\mathcal{O}(h^\alpha)$.
However, we observe that $h^{- \alpha} (q(h,t) - q^*(t))$ is $\mathcal{O}(1)$, suggesting that this quantity can be modelled using a stationary \ac{GP}.
We therefore take $E \sim \mathcal{GP}(0, \sigma^2 \rho_E k_E)$, where $\rho_E > 0$ is a parameter to be estimated, and 
\begin{equation}
	 k_E((h,t),(h',t')) = (h h')^\alpha \psi\left( |h-h'| / \ell_h \right) \cdot k_G(t,t') \label{eq: order cov}
\end{equation}
for a radial basis function $\psi$, scaled to satisfy $\psi(0) = 1$, and a length-scale parameter $\ell_h > 0$ to be estimated.
Note how \eqref{eq: order cov} separates the $h$ and $t$ dependence of $E$ in the prior, and adopts the same covariance model $k_G$ that was used to model the $t$ dependence of $G$. 
This can be motivated by the alternative perspective that follows from observing that $Q$ is a \ac{GP} with covariance function
\begin{align}
k_Q((h,t),(h',t')) = \sigma^2 \left\{ b(t) \cdot b(t') + \rho_G k_G(t,t') \left( 1 + \rho_E \frac{k_E((h,t),(h',t'))}{k_G(t,t')} \right) \right\} \label{eq: Q cov fn}
\end{align}
where $k_E/k_G$ is a kernel only depending on $h$. 
Written this way, the model is seen to perform universal kriging over $T$ with a covariance adjusted by a multiplicative error arising from non-zero values of $h$. 

\paragraph{Higher-order convergence:}
The \ac{GP} specification just described is not arbitrary; it ensures that the higher-order convergence property of \ac{RDAL} is realised in \ac{BBPN}.
Consider again the setting in \Cref{prop: richard}.
Suppose that there exist $L, \varepsilon_0 > 0$ and $\beta \in (0,1]$ such that $\lvert 1 - \psi(\varepsilon) \rvert \leq L \, \varepsilon^\beta$ for all $\varepsilon \in [0, \varepsilon_0)$.
Then the posterior mean $\mathbb{E}[Q(0)|D_h]$ satisfies $\lvert q^* - \mathbb{E}[Q(0)|D_h] \rvert = \mathcal{O}(h^{\alpha+\beta})$ as $h \to 0$.
Thus if $\psi$ is Lipschitz (\emph{i.e.} $\beta = 1$), \ac{BBPN} achieves the same higher-order convergence, $\alpha + 1$, as \ac{RDAL}.
In this context we recall that any Mat\'{e}rn covariance function of smoothness at least $1/2$ is Lipschitz.
The proof is provided in \Cref{ap: proof-error}.

\paragraph{Model parameters:} 
The free parameters of our prior model, $\sigma^2$, $\rho_G$, $\rho_E$, $\ell_h$, and the $\ell_{t,i}$ for $i = 1, \dots,p$, are collectively denoted $\theta$.
For all experiments in this article, $\theta$ was set using the maximum likelihood\footnote{Alternative approaches, such as cross-validation, could also be used; see Chapter 5 of \cite{rasmussen06}.
Our choice of maximum likelihood was motivated by the absence of any degrees of freedom (such as the number of folds of cross-validation), which permits a more objective empirical assessment. } estimator $\theta_{\textsc{ML}}$.
Our choice of parameterisation ensures that the maximum likelihood estimate for the overall scale $\sigma^2$, denoted $\sigma^2_{\textsc{ML}}$, has a closed form expression in terms of the remaining parameters. This is analytically derived in \Cref{app: ml explain} and can be plugged straight into the likelihood.
Gradients with respect to the remaining $3+p$ parameters are derived in \Cref{app: ml explain}, and gradient-based optimisation of the log-likelihood was implemented for the remaining parameters. 

\begin{remark} \label{rem: identify}
\ac{GP} interpolation, as with classical \ac{RDAL}, is not parameterisation invariant. 
Thus some care is required to employ a parameterisation of $h$ that is amenable to the construction of a \ac{GP} interpolant.
The effect of differing parameterisations is explored in \Cref{app: prior sensitivity}.
\end{remark}

\begin{remark}
The classical definition of \ac{RDAL} presupposes that, in order to employ the method, the order $\alpha$ must be known \emph{a priori} \cite{burg09}. 
However if $\alpha$ is not known, the probabilistic perspective affords us the opportunity to learn $\alpha$ as an additional parameter in the statistical model\textemdash a procedure with no classical analogue. 
The feasibility of learning $\alpha$ is explored in \Cref{subsec: eigenvalue}.
\end{remark}

\paragraph{Code:} Software for \ac{BBPN}, including code to reproduce the experiments in \Cref{sec: empirical}, can be downloaded from $\mathtt{github.com/oteym/bbpn}$.

\section{Experimental Assessment} \label{sec: empirical}
This section reports a rigorous experimental assessment of \ac{BBPN}.
Firstly, \Cref{subsec: odes} demonstrates that \ac{BBPN} is competitive with existing \ac{PN} methods in the context of \acp{ODE}.
This result is somewhat surprising, given the black box nature of \ac{BBPN} compared to the bespoke nature of existing \ac{PN} methods for \acp{ODE}.
Secondly, in \Cref{subsec: eigenvalue} we demonstrate the versatility of \ac{BBPN} by applying it to the nonlinear problem of eigenvalue computation, for which no \ac{PN} methods currently exist.
Finally, in \Cref{sec:PDEs} we use \ac{BBPN} to provide uncertainty quantification for state-of-the-art numerical methods that aim to approximate the solution of nonlinear \acp{PDE}.

\paragraph{Default Settings:}

We use Mat{\'e}rn(1/2) kernels for $\phi_i$ and $\psi$, \emph{i.e.} $\phi_i(t_i,t_i') = \exp (- \Vert t_i - t_i'\Vert / \ell_{t,i} )$,
and similarly \emph{mutatis mutandis} for $\psi$. 
These kernels impose a minimal continuity assumption on $q$ without additional levels of smoothness being assumed.
Sensitivity of results to the choice of kernel is investigated in \Cref{app: prior sensitivity}.

\paragraph{Performance Metrics:}
\ac{PN} is distinguished from traditional numerical analysis by its aim to provide probabilistic uncertainty quantification, but nevertheless approximation accuracy remains important.
To perform an assessment on these terms, we considered two orthogonal metrics.
Firstly, we compute the \emph{error} of the point estimate (mean), denoted $W \coloneqq \Vert \mathbb{E}[Q(0,\cdot)|D] - q^*(\cdot) \Vert$, where the norm is taken over $t \in T'$ where $T'$ is either $T$ itself or a set of representative elements from $T$.
Secondly, and most importantly from the point of view of \ac{PN}, we consider the \emph{surprise} 
$S \coloneqq \Vert \mathbb{C}[Q(0,\cdot)|D]^{-1/2} ( \mathbb{E} [Q(0,\cdot)|D)] - q^*(\cdot)  ) \Vert$, where $\mathbb{C}[Q(0,\cdot)|D]$ denotes the posterior covariance matrix. If the true quantity of interest $q^*$ was genuinely a sample from $Q(0,\cdot)|D$, then $S^2$ would follow a $\chi^2$ distribution with $|T'|$ degrees of freedom. This observation enables the \emph{calibration} of a \ac{PN} method to be assessed \citep{cockayne2021testing}. 
Both metrics naturally require an accurate approximation to $q^*$ to act as the ground truth, which is available using brute force computation in \Cref{subsec: odes,subsec: eigenvalue} but not in \Cref{sec:PDEs}.
The role of \Cref{sec:PDEs} is limited to demonstrating \ac{BBPN} on a problem class that is challenging even for state-of-the-art methods.

\begin{figure}[!t]
\centering
\hspace*{-0.5em}\includegraphics[width=1.025\textwidth]{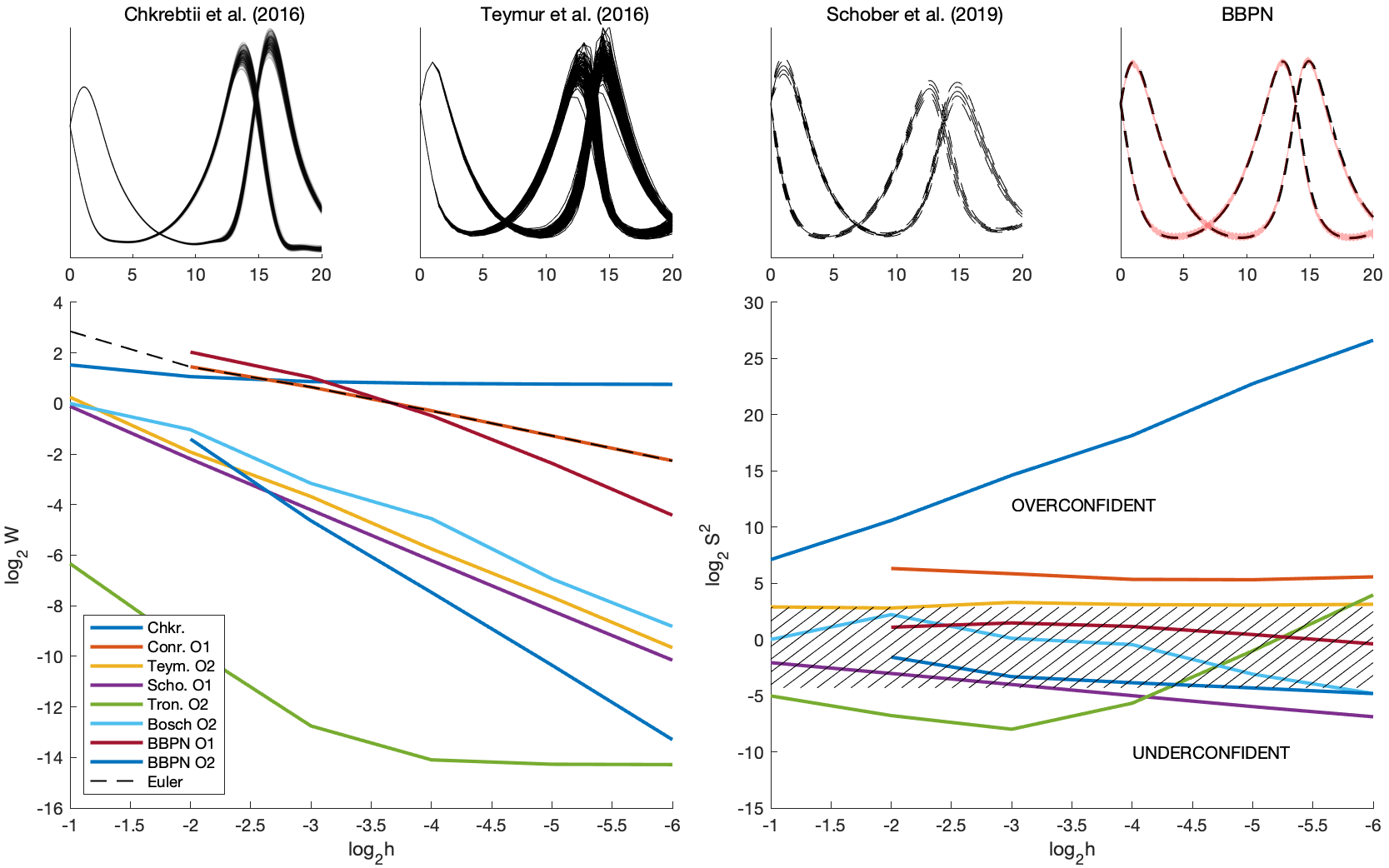}
\caption{Ordinary differential equations.
Top: Output from three existing \ac{PN} algorithms \cite{chkrebtii16,teymur16,schober18} and \ac{BBPN}, applied to the Lotka--Volterra \ac{IVP}.  
Bottom left: The error $\log_2 W$ at the final time point $t_\mathrm{end}=20$, as a function of the time step size $h$.
Bottom right: The surprise $\log_2 S$ at $t_\mathrm{end}=20$, with the central 95\% probability band of a $\chi_2^2$ random variable shaded.
Methods shown with (where applicable) their order: Chkr. \cite{chkrebtii16}; Conr. O1 \cite{conrad16}; Teym. O2 \cite{teymur16}; Scho. O1 \cite{schober18}; Tron. O2 \cite{tronarp18}; Bosch O2 \cite{bosch21}; BBPN O1 \& O2; and (traditional) Euler. }
\label{fig:ode figure}
\end{figure}

\subsection{Ordinary Differential Equations} \label{subsec: odes}

The numerical solution of \acp{ODE} has received considerable attention in \ac{PN}, with several sophisticated methods available to serve as benchmarks.
Here we consider numerical solution of the following Lotka--Volterra \ac{IVP}, a popular test case in the \ac{PN} literature:
\[
\frac{\mathrm{d}\mathbf{y}}{\mathrm{d}t} = f(t,\mathbf{y}) = \left[  \begin{array}{c} 0.5 y_1 - 0.05 y_1y_2 \\ - 0.5 y_2 + 0.05 y_1y_2 \end{array} \right],\ \quad \mathbf{y}(0) = \left[ \begin{array}{c} 20 \\ 20 \end{array} \right] 
\]
The aim in what follows is to approximate the quantity of interest $q^* = \mathbf{y}(t_{\text{end}})$ for $t_{\text{end}} = 20$.
The top row of \Cref{fig:ode figure} displays output from three distinct \ac{PN} methods due to \cite{chkrebtii16,teymur16,schober18}, as well as BBPN.
(For these plots the coarse step-size $h=0.5$ was used, so the probabilistic output can be easily visualised.)
In each case, these methods treat a sequence of evaluations of the gradient field $f$ as data which are used to constrain a random variable model for the unknown solution of the \ac{ODE}.
Their fundamentally differing character makes direct comparisons challenging, particularly if we are to account for computational cost. 
However, each algorithm has a recognisable discretisation parameter $h$, so it remains instructive to study their $h \rightarrow 0$ limit. (In most cases $h$ represents a time step size, but the method of \cite{bosch21} is step-size adaptive; in this case $h$ is an error tolerance that is user-specified.)
The methods of \cite{conrad16}, \cite{chkrebtii16}, and \cite{teymur16} require parallel simulations to produce empirical credible sets, and thus have a significant computational cost. 
The methods of \cite{schober18}, \cite{tronarp18} and \cite{bosch21} are based on Gaussian filtering and are less computationally demanding, though in disregarding nonlinearities the description of uncertainty they provide is not as rich. 
Interestingly, the output from \cite{chkrebtii16} becomes overconfident as $h \rightarrow 0$, with $S^2$ being incompatible with a $\chi_2^2$ random variable, while the output from \cite{schober18} becomes somewhat pessimistic in the same limit.
Aside from these two outputs, the other PN methods considered appear to be reasonably calibrated.

To illustrate \ac{BBPN}, our data consist of the final states produced by either an Euler (order 1) or an Adams--Bashforth (order 2) algorithm, which were performed at different resolutions $\{h_i = 2^{-i} , i = 1\dots,6\}$. The dataset\footnote{In this experiment the two components of $q^*$ were treated as \emph{a priori} independent, but this is not a specific requirement of \ac{BBPN} and dependence between outputs can in principle also be encoded into the \ac{GP} model.} is augmented cumulatively, so that for $i=i'$, all data generated by runs $1,\dots,i'$ are used. The finest resolution in each case, $h_i$, is simply denoted $h$.
For this experiment we use a prior with constant intercept, \emph{i.e.} $v=1$ and $b_1(t) = 1$. The \ac{BBPN} output, shown in the bottom row of \Cref{fig:ode figure}, is observed to be calibrated, 
and the (order 2) output provides the most accurate approximation among all calibrated PN methods considered. 
Note in particular how \ac{BBPN} accelerates the convergence of the Euler method from first order to second order, akin to \ac{RDAL}. 
In terms of computational cost, \ac{BBPN} requires running a traditional numerical method at different resolutions, as well as the fitting of a \ac{GP}. 
In this experiment, the computational cost of \ac{BBPN} was intermediate between the filtering approach of \cite{schober18} and the sampling approaches of \cite{chkrebtii16} and \cite{teymur16}. 
Further details, including the sources of all these codes, are given in \Cref{app: ODEs}.

\begin{figure}[!t]
\centering
\hspace*{-1em}\includegraphics[width=1.04\linewidth]{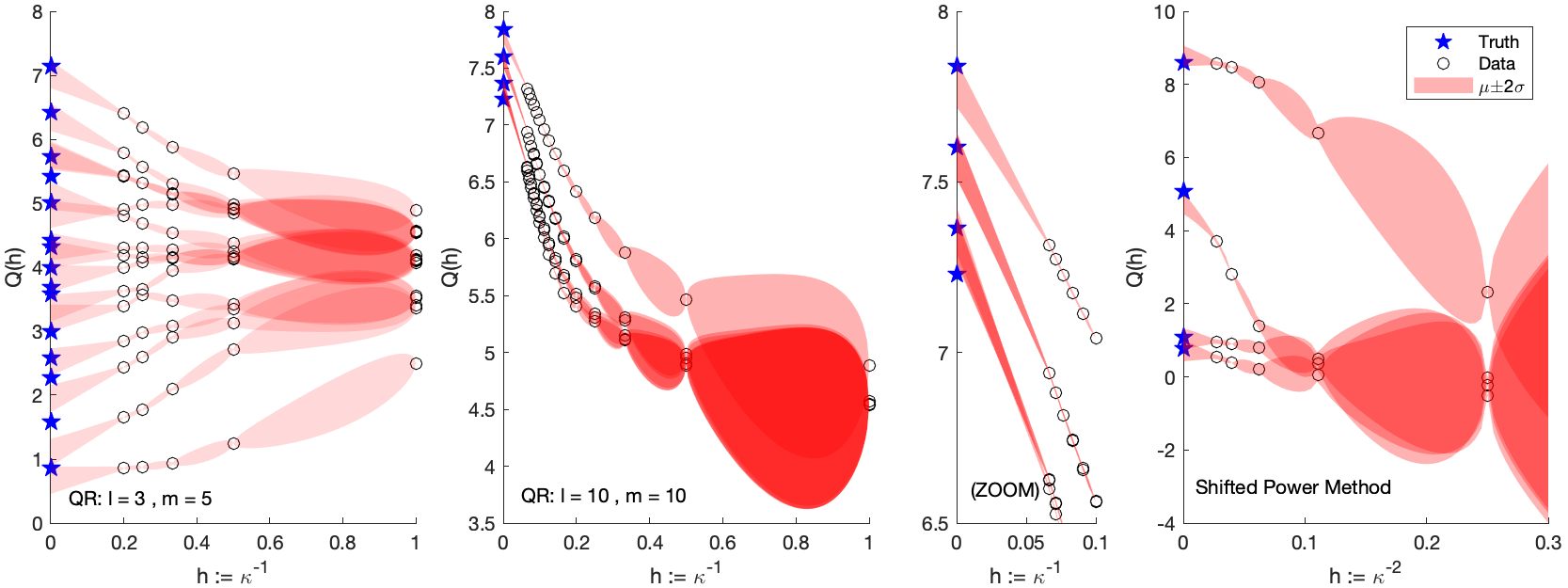}
\captionof{figure}[]{Eigenvalue problems. Left and centre: QR algorithm.
Right: Shifted power method. 
Details of each simulation are given in the main text. All plots show red shaded $\pm2\sigma$ credible intervals, numerical data as black circles, and true eigenvalues as blue stars.  }
\label{fig:eig figure}
\end{figure}

\subsection{Eigenvalue Problems} \label{subsec: eigenvalue}

The calculation of eigenvalues is an important numerical task that has not yet received attention in \ac{PN}. 
In this section we apply \ac{BBPN} to (1) the QR algorithm for matrix eigenvalue discovery, and (2) a cutting-edge adaptive algorithm for the \emph{tensor} eigenproblem, called the shifted power method \cite{kolda11}. 
In these examples, the order $\alpha$ is \emph{unknown} and we append it to $\theta$ as an additional parameter to be estimated using maximum likelihood.

\paragraph{QR Algorithm:} Take $T = \{1,\dots,n\}$. Let $w\in \mathbb{N} $ and define $h \coloneqq w^{-1} $.
Given a matrix $A \in  \mathbb{R}^{n\times n}$, its Schur form $A_\infty$ is approximated by matrices $A_w \coloneqq R_{w-1}Q_{w-1}$, where $R_{w-1}$ and $Q_{w-1}$ arise as a result of performing a QR decomposition on $A_{w-1} = Q_{w-1}R_{w-1}$, and where $A_0 \coloneqq A$.
Then $q(h,\cdot)$ is the vector $\mathrm{diag}(A_{h^{-1}})$, and $q^\ast =\mathrm{diag}(A_{\infty})$ is the vector of eigenvalues of $A$. 
As a test problem, whose eigenvalues are available in closed form (see \Cref{app: eigenvalues}), we consider the following family of sparse matrices that arise as the discrete Laplace operator in the solution of the Poisson equation by a finite difference method with a five-point stencil. Let the $l\times l$ matrix $B$ and the $ml\times ml$ matrix $A$ be defined by
\[
B =   \begin{pmatrix}
   4 &  -1 &  &  \\
   -1 &  4 &   -1  & \\
    & \ddots &  \ddots &  -1  \\
    &    &   -1  & \hphantom{-} 4 
    \end{pmatrix}\ , \qquad
    A = \begin{pmatrix}
   B &  -I &  &  \\
   -I &  B &   -I  & \\
    & \ddots &  \ddots &  -I  \\
    &    &   -I  & \hphantom{-} B 
    \end{pmatrix} .
\]

\ac{BBPN} output for this problem is displayed in \Cref{fig:eig figure}. In the left-hand pane, we take $l=5,m=2$ and perform 5 QR iterations, displaying all 10 eigenvalues. In the centre pane, we take $l=10,m=10$ and perform 15 iterations. For clarity, this pane only displays the largest few eigenvalues of this $100\times 100$ matrix, and we also show a zoomed-in crop to better demonstrate the extrapolation quality. Both examples show the convergence of $Q(h,\cdot)$ to $q^\ast(\cdot)$ as $w \rightarrow \infty$. Recall that $\alpha$ is \emph{inferred} in these simulations \textemdash the maximum likelihood values were, respectively, 1.0186 and 1.0167.

The extrapolation performed by our \ac{GP} model is seen visually to be effective and almost all true eigenvalues are contained within the $\pm 2\sigma$ credible intervals plotted. 
For comparison, in \Cref{app: prior sensitivity} we contrast the result of using a \emph{stationary} \ac{GP} model (\emph{i.e.} $\alpha=0$). The extrapolating properties of that \ac{GP} are immediately seen to be unsatisfactory, and we support this observation by examining the calibration of the two approaches, in a similar manner to in \Cref{subsec: odes}. This analysis strongly supports our proposed \ac{GP} specification in Section \ref{subsec: GPs}.

\paragraph{Shifted Power Method:} This iterative algorithm, due to \cite{kolda11} and implemented in \cite{bader21}, finds (random) eigenpairs of higher-order tensor systems, and we include it to demonstrate \ac{BBPN} on a challenging problem in linear algebra. 
For $\bm{\mathcal{A}}$ a symmetric $m$th-order $n$-dimensional real tensor, and $\bm{\mathrm{x}}$ an $n$-dimensional vector, define
\[
\left(\bm{\mathcal{A}}\bm{\mathrm{x}}^{m-1}\right)_{i_1} \coloneqq \sum\limits_{i_2 = 1}^n \cdots \sum\limits_{i_m = 1}^n a_{i_1i_2\dots i_m}x_{i_2}\dots x_{i_m}\ , \qquad i_1 = 1,\dots,n ,
\]
and say that $\lambda \in \mathbb{R}$ is an eigenvalue of $\bm{\mathcal{A}}$ if there exists $\bm{\mathrm{x}} \in \mathbb{R}^n$ such that $\bm{\mathcal{A}}\bm{\mathrm{x}}^{m-1} = \lambda \bm{\mathrm{x}}$ and $\bm{\mathrm{x}}^\top\bm{\mathrm{x}} = 1$. 
Here we take $n=m=6$ and produce a random symmetric tensor using the $\mathtt{create\_problem}$ function of \cite{bader21}. 
Two parameterisations of $q$ were considered; $h \coloneqq w^{-1}$ and $h \coloneqq w^{-2}$, where $w$ denotes the number of iterations performed, with results based on the latter parameterisation presented in the right-hand pane of \Cref{fig:eig figure}. (The maximum likelihood value for $\alpha$ in this example was 1.3318.)
It can be seen that, after 5 iterations, \ac{BBPN} is more accurate than each of the individual approximations on which it was trained.
The choice of parameterisation affects the performance of \ac{BBPN}, an issue we explore further in \Cref{app: prior sensitivity}. 

In this example there is no additional computational cost to \ac{BBPN} in the data collection stage, since the dataset is generated during a single run of an iterative numerical method. 
Therefore the only overhead is due to fitting the \ac{GP}; though for this example this cost is itself negligible.
All details, including a systematic assessment of error $W$ and surprise $S$ as $h$ is varied for both the QR algorithm and the shifted power method, are given in \Cref{app: eigenvalues}.

\begin{remark}
In this section we have implicitly modelled eigenvalues as \emph{a priori} independent, for simplicity of exposition. Heuristics from random matrix theory suggest that, when treated probabilistically, eigenvalues may be better modelled with some non-trivial dependence structure. We note that this additional structure can be easily incorporated into the prior GP model for BBPN.
\end{remark}

\subsection{Partial Differential Equations} \label{sec:PDEs}
To demonstrate the potential of \ac{BBPN} on a challenging problem for which state-of-the-art numerical methods are required, we consider numerical solution of the \emph{Kuramoto--Sivashinsky equation} \cite{Kuramoto78,Sivashinsky77}  
\begin{equation}\label{eq:KSE}
\partial_t u + \partial_x^4 u + \partial_x^2 u + u \partial_x u = 0.
\end{equation}

This equation is used to study a variety of reaction-diffusion systems, producing complex spatio-temporal dynamics which exhibit temporal chaos\textemdash the characteristics of which depend strongly on the amplitude of the initial data and the domain length.
We consider a flat initial condition $u(x,0) = \exp(-0.01x^2)$ with periodic boundary condition on the domain $0 \leq x \leq 1$, and numerically solve to $t = 200$. Our quantity of interest is therefore $u(x,200)$ for the domain $x \in [0,1]$.

\begin{figure}[!t]
    \centering
   \hspace*{-1em} \includegraphics[width=1.03\textwidth]{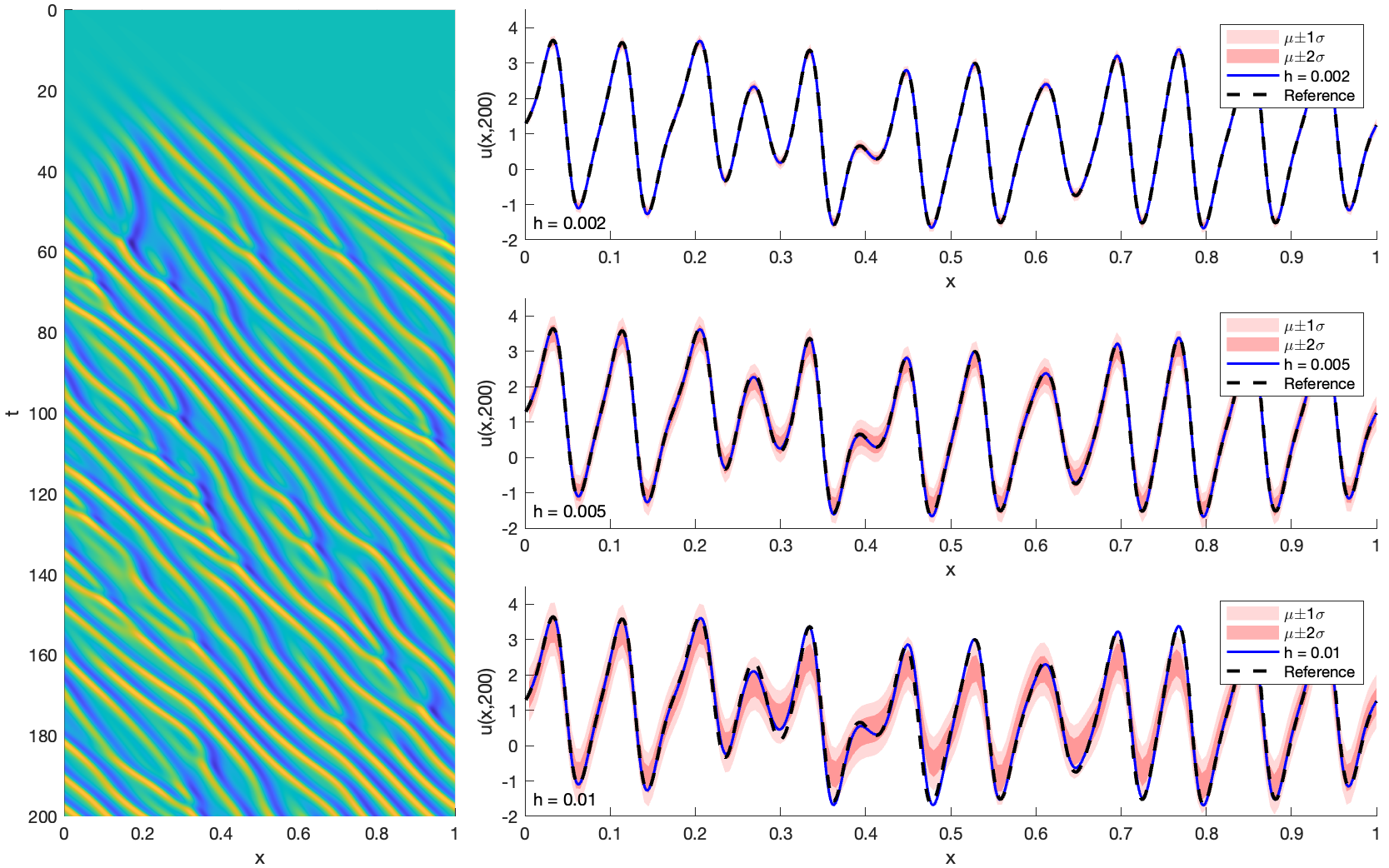}
    \caption{Partial differential equations.
    Left: Solution to the Kuramoto--Sivashinsky equation. 
    Right: Approximation of the solution at the final time point ($t=200$) using \ac{BBPN}, based on minimum time step sizes $h\in \{0.002, 0.005, 0.01\}$.
    Posterior mean (blue) and credible regions (shaded) are displayed.
    A reference solution (dashed black) is obtained by taking $h = 0.0005$.
    }
    \label{fig:KSE}
\end{figure}

To obtain numerical solutions to \eqref{eq:KSE} we transfer the problem into Fourier space and apply the popular fourth-order time-differencing ETD RK4 numerical scheme; see \Cref{app: PDEs} and \cite{Kassam2005}. 
ETD RK4 was designed to furnish fourth-order numerical solutions to time-dependent \acp{PDE} which combine low-order nonlinear terms with higher-order linear terms, as in \eqref{eq:KSE}. 
Computing accurate approximations to the solution of chaotic, stiff \acp{PDE} is a challenging problem for existing \ac{PN} methods because computationally demanding high-order approximations across both spatial and temporal domains are required. 
Here, we assess \ac{BBPN} applied to three sequences of five runs of ETD RK4, with minimum temporal step size $h = \delta t$ and, for simplicity, a fixed spatial step size $\delta x=0.001$ throughout.\footnote{For the $h=0.002$ simulation in \Cref{fig:KSE}, we have $h_i \in \{ 0.002,0.005,0.01,0.02,0.05\}$, for the $h=0.005$ simulation we have $h_i \in \{0.005,0.01,0.02,0.05,0.1\}$, and for the $h=0.01$ simulation we have $h_i \in \{0.01,0.02,0.05,0.1,0.2\}$}
A reference solution was generated by taking $h=0.0005$, but this cannot of course be guaranteed to be an accurate approximation to the true solution of \eqref{eq:KSE}. 

Results shown in \Cref{fig:KSE} are encouraging; not only can accurate approximations be produced, but the associated uncertainty regions appear to be reasonably well calibrated, insofar as the magnitude of the uncertainty is consistent with the magnitude of the discrepancy between the posterior mean and the reference solution. Full details of these simulations are contained in \Cref{app: PDEs}.

\section{Discussion} \label{sec: discuss}

This paper presented \textit{\acl{BBPN}}, a simple yet powerful framework that bridges the gap between existing \ac{PN} methods and the numerical state-of-the-art.
Positive results were presented on the important problems of numerically approximating \acp{ODE}, eigenvalues, and \acp{PDE}.
Our main technical contribution is a probabilistic generalisation of Richardson's deferred approach to the limit, which may be of independent interest.

The main drawbacks, compared to existing \ac{PN}, are a possibly increased computational cost and the additional requirement to model the error of a traditional numerical method. 
Compared to existing \ac{PN}, in which detailed modelling of the inner workings of numerical algorithms are exploited, only the order of the numerical method is used in \ac{BBPN} (and we can even dispense with that, as in \Cref{subsec: eigenvalue}), which may reduce its expressiveness in some settings. 
However, despite the black box approach, \ac{BBPN} was no less accurate than existing \ac{PN} in our experiments, and in fact the higher-order convergence property may enable \ac{BBPN} to out-perform existing \ac{PN}.

Some avenues for further research (that we did not consider in the present article) include the use of more flexible and/or computationally cheaper alternatives to \acp{GP}, the adoption of principles from experimental design to sequentially select resolutions $h_i$ given an overall computational budget, and the simultaneous use of different traditional (or even probabilistic) numerical methods within \ac{BBPN}.

\vfill

\begin{ack}
This work was supported by the Lloyd's Register Foundation programme on data-centric engineering at the Alan Turing Institute, UK.
The authors wish to thank Jon Cockayne and Ilse Ipsen for feedback on earlier versions of the manuscript, and Nicholas Kr\"{a}mer for assistance with the \texttt{probnum} package.
\end{ack}

\bibliographystyle{ieeetr}
\interlinepenalty=10000
\bibliography{bibliography}

\newpage
\appendix
\section*{Supplementary Material}
\raggedbottom
\setcounter{page}{1}
\pagenumbering{roman}

These appendices contain supplementary material for the paper \textit{Black Box Probabilistic Numerics}.

\section{Proof of \Cref{prop: richard}}
\label{ap: richard proof}

The equation of the straight line through two points $(x_1,y_1)$ and $(x_2,y_2)$ is given by 
\[
\frac{y-y_1}{x-x_1} = \frac{y_2 - y_1}{x_2 - x_1} .
\]
Substituting the points $(h^\alpha,q(h))$ and $((\gamma h)^\alpha,q(\gamma h))$, and taking $x=0$, we have
\[
y = q(h) - \frac{q(\gamma h) - q(h)}{\gamma^\alpha - 1} .
\]
By the assumption that $q$ is of order $\alpha$, we have the expansions $q(h) = q^\ast + Ch^\alpha + \mathcal{O}(h^{\alpha+1})$ and $q(\gamma h) = q^\ast + C(\gamma h)^\alpha + \mathcal{O}(h^{\alpha+1})$, and then by substitution and straightforward cancellation we find
\[
y = q^\ast + \mathcal{O}(h^{\alpha+1}) .
\]
Therefore the $y$-intercept of the line is an approximation of $q^\ast$ of order $\alpha+1$.

\section{Gaussian Processes for BBPN}
\label{ap: learning}

This appendix contains full details of how analytic conditioning formulae are obtained and how maximum likelihood estimates are calculated.

\subsection{Conditioning Formulae}
\label{app: conditioning formulae}

It will be convenient to introduce \textit{lexicographic ordering}, where the indices 
\begin{equation}
\{ (i,j) : j = 1,\dots,m_i, \; i = 1,\dots,n \} \label{eq: indices}
\end{equation}
are ordered first by $i$ and then, for indices with the same $i$, by $j$.
Let $h_{(l)}$ and $t_{(l)}$ denote, respectively, the values of $h_i$ and $t_{i,j}$ corresponding to the $l$'th ordered pair $(i,j)$ in \eqref{eq: indices}. 
Let $\mathbf{q}$ represent a column vector of length $m := \sum_{i=1}^n m_i$, with entries $\mathbf{q}_{\;\!l} := q(h_{(l)},t_{(l)})$ in lexicographic order.

From \eqref{eq: Q cov fn}, the prior model for $Q$ described in \Cref{subsec: GPs} has covariance function
\begin{equation}
    k_Q((h,t),(h',t')) = \sigma^2 [ b(t) \cdot b(t') + \rho_G k_G(t,t')+ \rho_E k_E((h,t),(h',t'))  ] , \label{eq:combined process}
\end{equation}
where the additivity follows from the assumptions $Q^* \indep E$ and $Z \indep G$.
Let $K_Q$ be an $m \times m$ matrix and $\mathbf{k}_Q(h,t)$ be an $m \times 1$ column vector with entries of the form
\begin{align}
(K_Q)_{l,l'} := k_Q((h_{(l)},t_{(l)}),(h_{(l')},t_{(l')})) \ ,\quad 
(\mathbf{k}_Q(h,t))_l := k_Q((h_{(l)},t_{(l)}),(h,t)) . \label{eq: def matrices}
\end{align}
Then standard Gaussian conditioning formulae (eg. Equation 2.19 in \cite{rasmussen06}) demonstrate that the conditional process $Q | D$ has mean and covariance functions
\begin{align}
    \mu_{Q | D}(h,t) & = \mathbf{k}_Q(h,t)^\top K_Q^{-1} \mathbf{q} \label{eq: post mean} \\
    k_{Q | D}((h,t),(h',t')) & = k_Q((h,t),(h',t')) - \mathbf{k}_Q(h,t)^\top K_Q^{-1} \mathbf{k}_Q(h',t')  \label{eq: post cov}
\end{align}
The mean and covariance functions of the marginal process $Q(0,\cdot) | D$ are extracted by setting $h$ equal to 0 in \Cref{eq: post mean,eq: post cov}.

\subsection{Proof of Higher-Order Convergence Result in \Cref{subsec: GPs}} \label{ap: proof-error}

For a scalar quantity of interest, the full covariance function in~\eqref{eq: Q cov fn} is
  \begin{equation*}
    k_Q(h, h') = a_1 + a_2 (h h')^\alpha \psi\bigg( \frac{\lvert h - h' \rvert}{\ell_h} \bigg)
  \end{equation*}
  for certain positive constants $a_1$ and $a_2$.
  For $\gamma \in [0,1]$, denote
  \begin{equation*}
    \psi_{h} = \psi\bigg( \frac{(1 - \gamma)h}{\ell_h} \bigg).
  \end{equation*}
  Then the conditional mean at $h=0$, given the data $D_h = \{ (h, q(h)), (\gamma h, q(\gamma h))\}$, is
  \begin{equation*}
    \begin{split}
      \mathbb{E}[Q(0)|D_h] &= \begin{pmatrix} q(h) \\ q(\gamma h) \end{pmatrix}^\top \begin{pmatrix} a_1 + a_2 h^{2\alpha} & a_1 + a_2 \gamma^\alpha \psi_h h^{2\alpha} \\ a_1 + a_2 \gamma^\alpha \psi_h h^{2\alpha} & a_1 + a_2 \gamma^{2\alpha} h^{2\alpha} \end{pmatrix}^{-1} \begin{pmatrix} a_1 \\ a_1 \end{pmatrix} \\
      &=  \frac{q(h) \gamma^\alpha(\gamma^\alpha  - \psi_h) + q(\gamma h) (1 - \gamma^\alpha \psi_h)}{ a_1 a_2 ( 1 - 2 \gamma^\alpha \psi_h + \gamma^{2\alpha} ) h^{2\alpha} + a_2^2 \gamma^{2\alpha} (1 - \psi_h^2) h^{4\alpha} } a_1 a_2 h^{2\alpha} \\
      &= \frac{q(h) \gamma^\alpha(\gamma^\alpha  - \psi_h) + q(\gamma h) (1 - \gamma^\alpha \psi_h)}{ a_1 ( 1 - 2 \gamma^\alpha \psi_h + \gamma^{2\alpha} ) + a_2 \gamma^{2\alpha} (1 - \psi_h^2) h^{2\alpha} } a_1.
      \end{split}
  \end{equation*}
  Inserting $q(h) = q^* + C h^\alpha + \mathcal{O}(h^{\alpha+1})$ and $q(\gamma h) = q^* + C \gamma^\alpha h^\alpha + \mathcal{O}(h^{\alpha+1})$ in the above equation yields
  \begin{equation*}
    \begin{split}
      \big| q^* - \mathbb{E}[Q(0)|D_h] \big| ={}& q^*\Bigg\lvert 1 - \frac{\gamma^\alpha(\gamma^\alpha  - \psi_h) + 1 - \gamma^\alpha \psi_h}{ a_1 ( 1 - 2 \gamma^\alpha \psi_h + \gamma^{2\alpha} ) + a_2 \gamma^{2\alpha} (1 - \psi_h^2) h^{2\alpha} } a_1 \Bigg\rvert \\
      &\qquad+ \Bigg\lvert \frac{\gamma^\alpha(\gamma^\alpha  - \psi_h) + \gamma^\alpha (1 - \gamma^\alpha \psi_h)}{ a_1 ( 1 - 2 \gamma^\alpha \psi_h + \gamma^{2\alpha} ) + a_2 \gamma^{2\alpha} (1 - \psi_h^2) h^{2\alpha} } \Bigg\rvert \lvert C \rvert \, a_1 h^{\alpha} \\
      &\qquad+ \Bigg\lvert \frac{\gamma^\alpha(\gamma^\alpha  - \psi_h) + 1 - \gamma^\alpha \psi_h}{ a_1 ( 1 - 2 \gamma^\alpha \psi_h + \gamma^{2\alpha} ) + a_2 \gamma^{2\alpha} (1 - \psi_h^2) h^{2\alpha} } \Bigg\rvert a_1 \mathcal{O}(h^{\alpha+1}) \\
      \leq{}& q^*\Bigg\lvert \frac{a_2 \gamma^{2\alpha} (1 - \psi_h^2) }{ a_1 ( 1 - 2 \gamma^\alpha \psi_H + \gamma^{2\alpha} ) + a_2 \gamma^{2\alpha} (1 - \psi_h^2) h^{2\alpha} } \Bigg\rvert h^{2\alpha} \\
      &\qquad+ \Bigg\lvert \frac{\gamma^\alpha (1 + \gamma^\alpha)}{ 1 - 2 \gamma^\alpha \psi_h + \gamma^{2\alpha} } \Bigg\rvert \lvert C \rvert \, \lvert 1 - \psi_h \rvert h^{\alpha} \\
      &\qquad + \mathcal{O}(h^{\alpha+1}).
      \end{split}
  \end{equation*}
  It follows from the Hölder assumption $\lvert 1 - \psi(\varepsilon) \rvert \leq L \, \varepsilon^\beta$ that $\lvert 1-\psi_h \rvert = \mathcal{O}(h^\beta)$.
  Therefore the second term, which dominates the right-hand side, is of order $\mathcal{O}(h^{\alpha+\beta})$.
  This concludes the proof.

\subsection{Maximum Likelihood Estimation}
\label{app: ml explain}

The parameters $\theta$ of the covariance function $k_Q$ are estimated from data using maximum likelihood. 
Recall that (with $\alpha$ known) $\theta$ consists of the parameters $\sigma$, $\rho_G$, $\rho_E$, $\ell_h$, and the $\ell_{t,i}$ for $i = 1, \dots,p$.
This parameterisation is deliberately chosen to enable the maximum likelihood estimator $\sigma_{\textsc{ML}}$ to be computed as an explicit function of the remaining components of $\theta$.
It is convenient to express
\begin{align*}
    k_Q((h,t),(h',t')) = \sigma^2 \overline{k}_Q((h,t),(h',t'))
\end{align*}
where $\overline{k}_Q((h,t),(h',t'))$ is \eqref{eq:combined process} with $\sigma = 1$.
Analogously define $\overline{K}_Q$ as in \eqref{eq: def matrices} but with $\sigma = 1$.
The log-likelihood of observing the dataset $D$ in \eqref{eq: triples of data} under the model for $Q$ defined in \eqref{eq: high level model} can then be expressed as
\begin{equation}
\mathcal{L}(\theta)  =  -\frac{m}{2}\log (2\pi) - m \log\sigma -\frac{1}{2}\log |\overline{K}_Q| - \frac{1}{2\sigma^2} \mathbf{q}^\top \overline{K}_Q^{-1} \mathbf{q} , \label{eq:log lik}
\end{equation}
where we note that $\overline{K}_Q$ does not depend on $\sigma$ but can depend on all the other components of $\theta$.
In the case of the overall amplitude parameter $\sigma$, it is possible to obtain an analytic expression for the value $\sigma_{\textsc{ML}}$ by differentiating and setting $\partial \mathcal{L}/\partial \sigma = 0$ \cite{karvonen20}. This gives
\begin{equation}
\sigma_{\textsc{ML}}^2 = \frac{ \mathbf{q}^\top \overline{K}_Q^{-1} \mathbf{q} }{m}
\end{equation}
Plugging $\sigma = \sigma_{\textsc{ML}}$ into \eqref{eq:log lik} gives
\begin{equation} 
\mathcal{L}(\theta | \sigma = \sigma_{\textsc{ML}}) =
- \frac{m}{2}  \log(\mathbf{q}^\top \overline{K}_Q^{-1} \mathbf{q}) - \frac{1}{2} \log|\overline{K}_Q| + C \label{eq: plug in likelihood}
\end{equation}
where $C$ is a constant in $\theta$.
From here, we employ numerical optimisation to maximise \eqref{eq: plug in likelihood} over the remaining $3+p$ degrees of freedom in $\theta$. 

It is important to ensure that numerical optimisation is successful, otherwise conclusions provided by BBPN could be an artefact of failure of the numerical optimisation method.
To this end, we undertake robust gradient-based optimisation on \eqref{eq: plug in likelihood}, using \texttt{MATLAB}'s packaged \texttt{fmincon} routine.
This requires calculation of the gradients of \eqref{eq: plug in likelihood} and explicit formulae will now be provided.

By differentiating \eqref{eq: plug in likelihood} we have 
\begin{equation}
\partial_\theta\mathcal{L}(\theta | \sigma = \sigma_{\textsc{ML}}) =  \frac{m}{2} \frac{ \mathbf{q}^\top \overline{K}_Q^{-1} (\partial_\theta \overline{K}_Q) \overline{K}_Q^{-1} \mathbf{q}}{ \mathbf{q}^\top \overline{K}_Q^{-1} \mathbf{q}} - \frac{1}{2} \text{tr}\big( \overline{K}_Q^{-1} (\partial_\theta \overline{K}_Q) \big)
\end{equation}
Define the matrices 
\begin{align*}
	(B)_{l,l'} & := b(t_{(l)}) \cdot b(t_{(l')})\ ,\ 
	(K_G)_{l,l'} := k_G(t_{(l)},t_{(l')})\ ,\ 
	(K_E)_{l,l'} := k_E((h_{(l)},t_{(l)}),(h_{(l')},t_{(l')})) .
\end{align*}
Then $\overline{K}_Q = B + \rho_G K_G + \rho_E K_E$, and it follows that
\begin{align*}
    \partial_{\rho_G} K_Q & =  K_G\ , & 
    \partial_{\ell_h} K_Q & =  \rho_E \partial_{\ell_h} K_E\ , \\
    \partial_{\rho_E} K_Q & =  K_E\ , & 
    \partial_{\ell_{t,i}} K_Q & =  \rho_G \partial_{\ell_{t,i}} K_G + \rho_E \partial_{\ell_t} K_E
\end{align*}
The low-level terms such as $\partial_{\ell_h} K_E$ can readily be computed by hand and will depend on the radial basis functions $\phi_i$ and $\psi$ adopted in $K_G$ and $K_E$.
Note that if $\alpha > 0$ is treated as unknown and appended to the parameter vector $\theta$, as in \Cref{subsec: eigenvalue}, a similar calculation can be performed to obtain the gradient with respect to $\alpha$ of \eqref{eq: plug in likelihood}.

The convergence of this gradient-based optimisation approach to a minimum of $\mathcal{L}(\theta)$ is verified empirically in \Cref{subsec: param ide}.

\begin{figure}[t]
    \centering
    \hspace*{-1em}\includegraphics[width=1.03\linewidth]{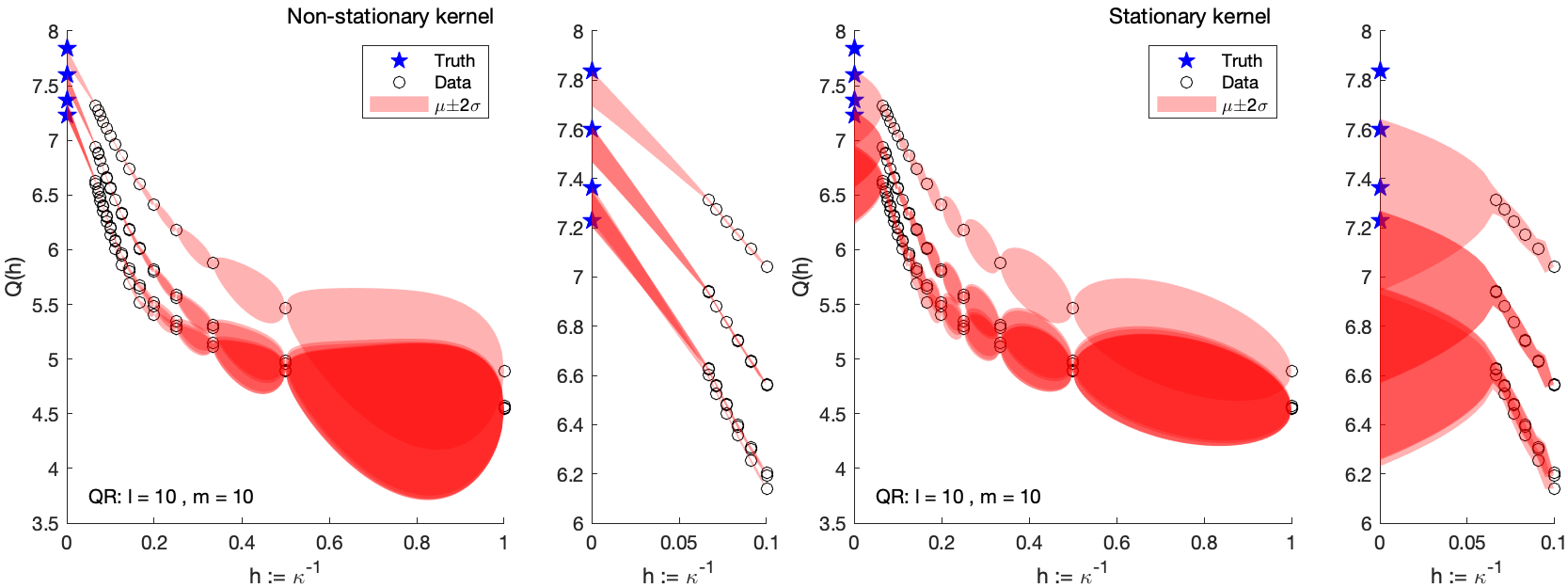}
    \caption{Comparison of stationary \eqref{eq: order cov}, on left, and non-stationary \eqref{eq: non-sta kernel}, on right, covariance functions for the QR algorithm example detailed in \Cref{subsec: eigenvalue}.}
    \label{fig:stationary}
\end{figure}

\section{Details of Empirical Assessment} \label{app: code sources}

This appendix contains full details for all experiments described in the main text.

\subsection{Riemann Sum Illustration in \Cref{fig:classical example}}
\label{app: Riemann}
\Cref{fig:classical example} considers the function 
$
f(x) = \sin^2(4\pi x) + \exp(x) - \tfrac{5}{2}x^4 + \tfrac{1}{2}\cos(16\pi x) + \tfrac{1}{4}\cos(20\pi x)
$.
The quantity of interest $q^\ast$ is the integral $\int_0^1 f(x) \d x$, which has the exact value $(\mathrm{e} - 1) \approx 1.71828$.

\ac{BBPN} was applied to the method of Riemann sums. 
The convergence of this method is first order, and we set $\alpha=1$ accordingly. We choose a range of step-sizes $h$ between $0.01$ an $0.08$, with the Riemann sum approximations plotted in the left pane of \Cref{fig:classical example}.
Hyperparameters of the \ac{GP} were set using maximum likelihood approach, as described in \Cref{app: ml explain}.

\subsection{Sensitivity to Prior Specification} \label{app: prior sensitivity}

In this section we consider the effect of varying several of the choices made during the specification of our prior model. 
The suitability of our non-stationary \ac{GP} model is considered in \Cref{app: stationarity}.
The effect of the choice of parametrisation for $h$ is considered in \Cref{app: stationarity}.
The choice of the kernel functions $\phi_i$ and $\psi$ is discussed in \Cref{app: choice of basis}.
Finally, the nature and number of the finite-dimensional basis terms $b_i$ is discussed in \Cref{app: choice of v}.
In each case we explore the impact of these aspects of the prior specification by reproducing figures from the main text under different settings within the \ac{GP} model.

\subsubsection{Stationary / Non-Stationary Error Model}
\label{app: stationarity}

Since the error $E(h,t)$ is assumed to vanish in the limit $h\rightarrow 0$, and since its scale is assumed to depend on the order $\alpha$ of the underlying numerical method, we specified a non-stationary \ac{GP} in \eqref{eq: order cov}. For the QR algorithm example in \Cref{subsec: eigenvalue}, we now contrast this with the same analysis performed with the stationary \ac{GP} whose covariance function is
\begin{equation}
	 \tilde{k}_E((h,t),(h',t')) = \psi\left( |h-h'| / \ell_h \right) \cdot k_G(t,t') \label{eq: non-sta kernel}
\end{equation}
\emph{i.e.} setting $\alpha = 0$ in \eqref{eq: order cov}.

From \Cref{fig:stationary} (right), we see that the extrapolation is extremely poor when a stationary \ac{GP} is used.
Moreover, the use of a stationary \ac{GP} leads in this case to over-confident predictions, with the true eigenvalues belonging outside of the $\pm 2\sigma$ credible intervals.
This provides strong support for the use of the non-stationary \ac{GP} that we propose in the main text.

\begin{figure}[t]
    \centering
    \hspace*{-1em}\includegraphics[width=1.03\linewidth]{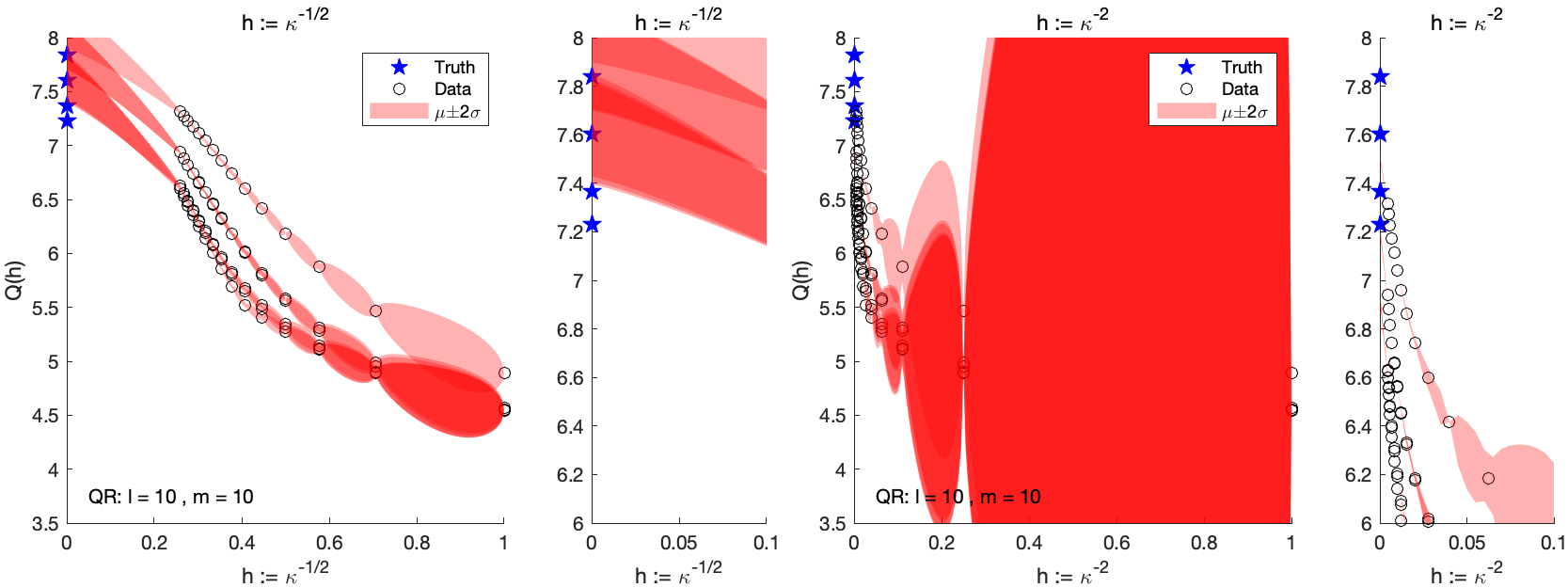}
    \caption{Comparison of different parameterisations for $h$ relative to the number of iterations $\kappa$ of the QR algorithm; $h := \kappa^{-1/2}$ (left); $h:= \kappa^{-2}$ (right)}
    \label{fig:alphachange}
\end{figure}

\subsubsection{Parameterisation of $h$}
\label{app: stationarity}

The choice of parameterisation of $h$ is also crucial to the operation of BBPN. While it is sometimes the case that an `obvious' parameterisation exists (such as the step-size in a time-stepping method, where the order $\alpha$ specifically refers to this quantity; or the overall tolerance level of a numerical method) this is, unfortunately, not always true. If some heuristic reasoning for determining this parameterisation is not available, we recommend some prior experimentation and comparison with calibration metrics such as surprise, introduced in \Cref{sec: empirical}.

For the QR algorithm example in \Cref{subsec: eigenvalue}, \Cref{fig:alphachange} shows the effect of replacing the parameterisation $h := \kappa^{-1}$ (as in \Cref{fig:eig figure}) with $h := \kappa^{-1/2}$ and $h := \kappa^{-2}$.
Although \ac{BBPN} continues to work, to an extent, with these alternative parametrisations, its predictive performance is somewhat diminished.

\subsubsection{Choice of Radial Basis Functions $\phi_i$ and $\psi$}
\label{app: choice of basis}

For all simulations in this article we specified Mat\'ern 1/2 kernels for $\phi_i$ and $\psi$. The motivation for this, stated in the preliminary notes in \Cref{sec: empirical}, is to impose the minimal continuity assumption on $q$ but not to assume additional levels of smoothness where this cannot be justified \emph{a priori}.

\begin{figure}[t]
    \centering
    \hspace*{-0.5em}\includegraphics[width=1.02\linewidth]{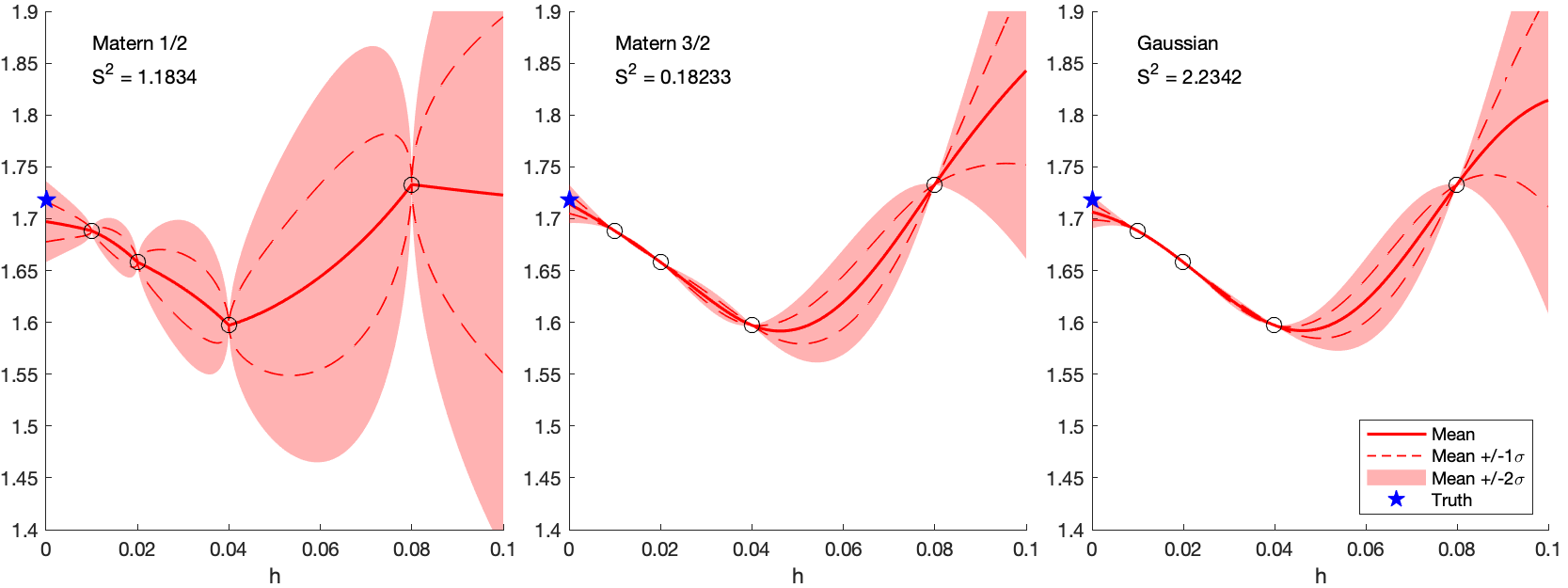}
    \caption{Comparison of different kernel types for the radial basis functions $\phi_i$ and $\psi$. Mat\'ern 1/2 (left); Mate\'ern 3/2 (centre); and Gaussian (right).}
    \label{fig:kernelchange}
\end{figure}

\Cref{fig:kernelchange} shows the effect of specifying instead Mat\'ern 3/2 or Gaussian kernels for $\phi_i$ and $\psi$ in the Riemann sum test problem in \Cref{fig:classical example}, contrasting with the Mat\'ern 1/2 kernel used there. In all cases, the same process of gradient-based optimisation was employed to automate the setting of the kernel hyperparameters. The additional smoothness of the mean interpolant is clearly visible in the higher Mat\'ern and Gaussian cases, but note also the difference in scale of the $\pm 2\sigma$ region.
In particular, the use of smoother kernels is associated with higher confidence in the predictive output, with the Gaussian kernel producing the largest value of the surprise $S^2$ (though this was still within the central $95\%$ region for a $\chi^2$ distribution, so we do not reject the hypothesis that the \ac{BBPN} output is calibrated).
On balance we err on the side of caution and recommend the Mat\'{e}rn 1/2 kernel for applications of \ac{BBPN}.

\subsubsection{Choice of Basis Functions $b_i$} 
\label{app: choice of v}

\begin{figure}[t]
    \centering
    \hspace*{-1em}\includegraphics[width=1.03\linewidth]{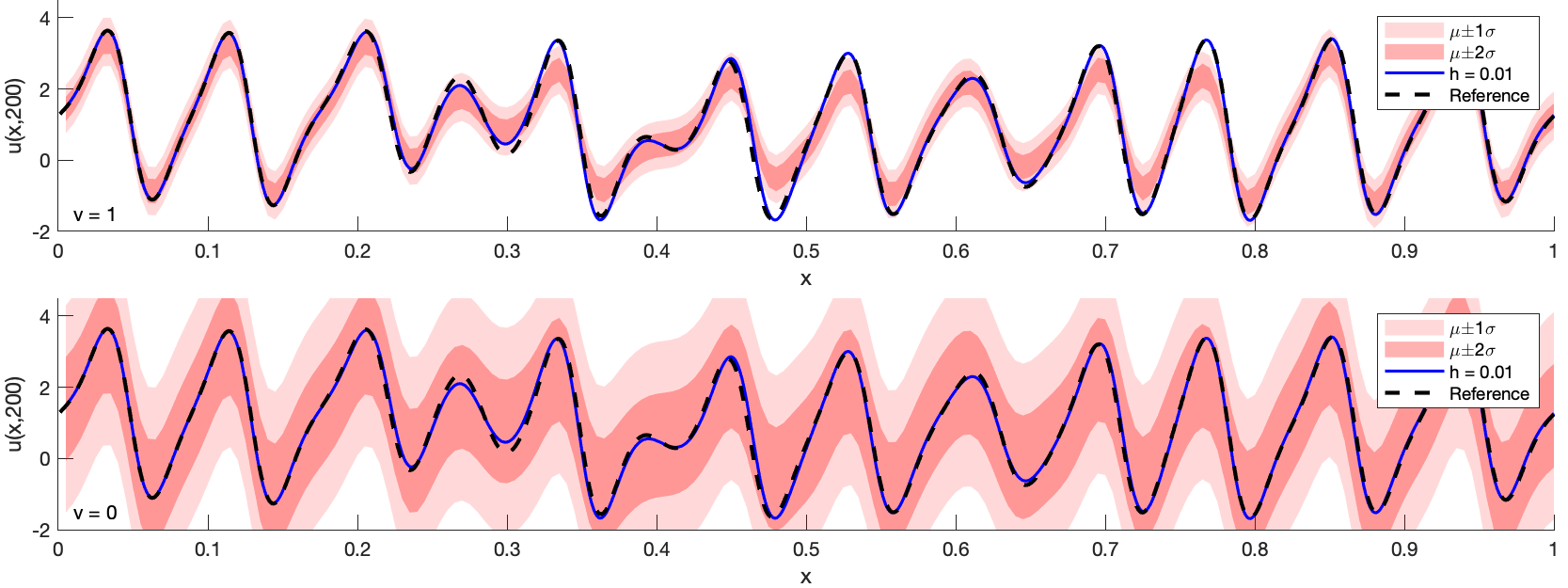}
    \caption{Comparison of the inclusion and exclusion of the first polynomial basis function (top: $v = 1$, bottom: $v= 0$) for the model in \Cref{sec:PDEs}.}
    \label{fig:vchange}
\end{figure}

In this section we demonstrate the purpose of including basis functions $b_i$ in the model for $G(t)$. To do so, we plot the output of the BPPN procedure for the PDE example in \Cref{fig:KSE}, since this example has non-trivial `$t$' domain (though the variable called $t$ in the model definition in \Cref{sec:setting} is in fact called $x$ here). The effect of including a constant basis function (\emph{i.e.} $v=1$ and $b_1(t) = 1$) is to allow the model a non-zero mean in $t$. For this example, the dynamics are mostly above the $0$ level and even a simple global mean would be more likely between $1$ and $2$. Omitting the basis function (\emph{i.e.} $v=0$), as shown in the bottom pane of \Cref{fig:vchange}, inflates the covariance to compensate for this misfit, and in this case results in an underconfident model.

In this example, it is unlikely that the additional inclusion of higher-order polynomial basis functions would be of use. Indeed our experiments showed this. However the oscilliating nature of the dynamics across the range of $t$ suggests a Fourier basis may be an appropriate mean model. Ideas along these lines are partially explored in \cite{kersting20a}, and a fuller investigation in the context of BBPN will be the subject of future work.

\subsection{Ordinary Differential Equations}
\label{app: ODEs}

Here we provide full details for the \ac{ODE} experiment in the main text.
In \Cref{subsec: ode detail} we explain how all the probabilistic \ac{ODE} solvers that we considered in the main text were implemented.
Then, in \Cref{subsec: param ide}, we present evidence that the gradient-based optimisation approach we employed to estimate the \ac{GP} hyperparameters in \ac{BBPN} has successfully converged.

\subsubsection{Details of Implementation} 
\label{subsec: ode detail}

In this section we describe in detail the sources and licences of the codes, as well as the settings used, to perform the comparison experiments in Section \ref{subsec: odes}. These codes are from different sources, span several years in release date, and are coded in different languages. They also accept inputs and give outputs in mutually inconsistent forms. This makes a `cloned-repository' solution from which results could be reproduced automatically impractical. In the interests of maximum possible transparency we manually collect and present code sources and parameter values here in the hope that interested readers will not find it difficult to reproduce our results locally if required. Recall that our simulations consist of varying input $h$.

The one-step-ahead sampling model of Chkrebtii \textit{et al.} \cite{chkrebtii16} (labelled `Chkr.' in \Cref{fig:ode figure}) was run using \texttt{MATLAB} code from \url{https://git.io/J33lL} with $\mathtt{nsolves}=100$, $\mathtt{N}=\lceil 20/h \rceil$, $\mathtt{nevalpoints}=1001$ and the $\mathtt{lambda}$ and $\mathtt{alpha}$ hyperparameters left at their default values (which depend on $\mathbb{N}$, and therefore $h$). This software has no explicitly-stated licence.

The perturbed integrator approach of Conrad et al. \cite{conrad16} and Teymur et al. \cite{teymur16} (labelled `Conr. O1' and `Teym. O2' was run using \texttt{MATLAB} code provided to us by the authors of the latter paper and not, as far as we are aware, publicly released.

The Gaussian filtering approach of Schober et al. \cite{schober18}, Tronarp et al. \cite{tronarp18} and Bosch et al. \cite{bosch21} (labelled `Scho. O1', `Tron. O2' and `Bosch O2') was run by installing the \texttt{Python} package \texttt{probnum} and using the function \texttt{probsolve\_ivp}. `Scho. O1' uses non-adaptive step-sizes and takes \texttt{algo\_order}$\ = 1$, and $\mathtt{method} = \mathtt{EK0}$; `Tron. O2' uses non-adaptive step-sizes and takes \texttt{algo\_order}$\ = 2$, and $\mathtt{method} = \mathtt{EK1}$; while `Bosch O2' uses \emph{adaptive} step-sizes and takes \texttt{algo\_order}$\ = 2$, and $\mathtt{method} = \mathtt{EK1}$. In the latter case, $h$ is taken as the relative tolerance $\mathtt{rtol}$ instead of the step-size. This software is Copyright of the ProbNum Development Team and is released under an MIT licence.

The reference solution used in calculating errors was calculated using \texttt{MATLAB}'s in-build \texttt{ode45} function with tolerances set using \texttt{odeset(`RelTol',3e-14,`AbsTol',1e-20)}

It is difficult to fairly compare the wall-clock times of these codes, particularly since they are written in different languages and are therefore run in different environments. For the example simulation in \Cref{fig:ode figure}, none of the examples took more than a few seconds on a 2018 MacBook Pro, and some were virtually instant. 
All publicly-available codes were downloaded or cloned on 22 April 2021.

\subsubsection{Parameter Identifiability}
\label{subsec: param ide}

\begin{figure}[t]
    \centering
    \hspace*{-1em}\includegraphics[width=1.03\linewidth]{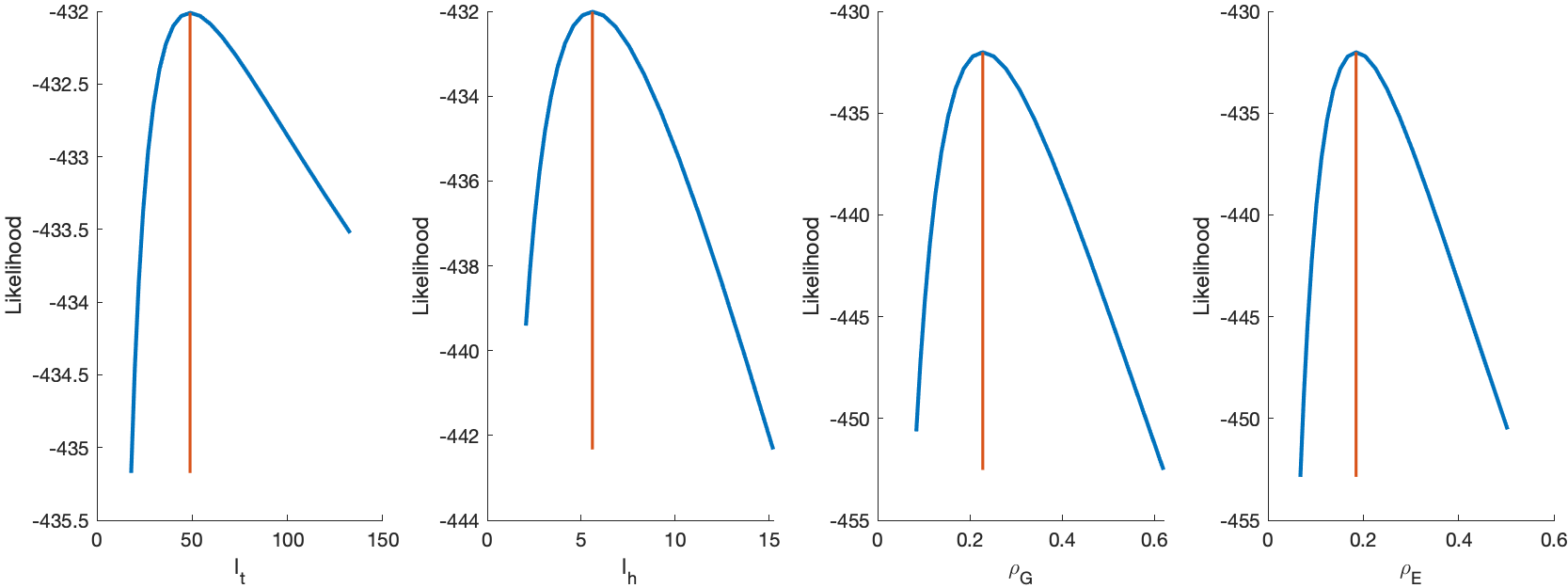}
    \caption{Likelihood variation in the neighbourhoods of the maximum likelihood values found by \texttt{MATLAB}'s \texttt{ode45} optimiser. In each case, the remaining parameters were fixed at their maximum likelihood value. The values determined by the optimiser are shown with a vertical orange line.}
    \label{fig:param ident}
\end{figure}

In order to assess the robustness of our gradient-based optimisation procedure for maximum likelihood estimation, we consider again the Lotka--Volterra model.
Here we will vary each parameter $\ell_t$, $\ell_h$, $\rho_G$ and $\rho_E$ in turn, holding all other parameters fixed at the values produced by the gradient-based optimisation method. 
The resulting plots are given in \Cref{fig:param ident}.

In this application, at least, we can be reasonably confident that the optimisation procedure has located a global maximiser in 4D (though, strictly speaking, we cannot confirm this from the univariate plots in \Cref{fig:param ident}). 
In general, and as is common in \ac{GP} modelling, model fit should always be assessed, in order to be confident in the data-driven nature of the \ac{GP} output, something particularly salient in numerical applications where calibration is of paramount importance.

\subsection{Eigenvalue Problems}
\label{app: eigenvalues}

In this section we provide certain further details for the eigenvalue problem presented in \Cref{subsec: eigenvalue}. We first note that the matrix $A$ defined there can be shown to have exact eigenvalues $4 - 2\cos (p\pi/(l+1)) - 2\cos (q\pi/(m+1)); p = 1,\dots,l ; q = 1,\dots,m$. The knowledge of the true values is required to facilitate the following analysis. For this section we take $l=3$ and $m=5$, as in the left-hand panes of \Cref{fig:eig figure}.

\begin{figure}[t]
    \centering
    \hspace*{-1em}\includegraphics[width=1.03\linewidth]{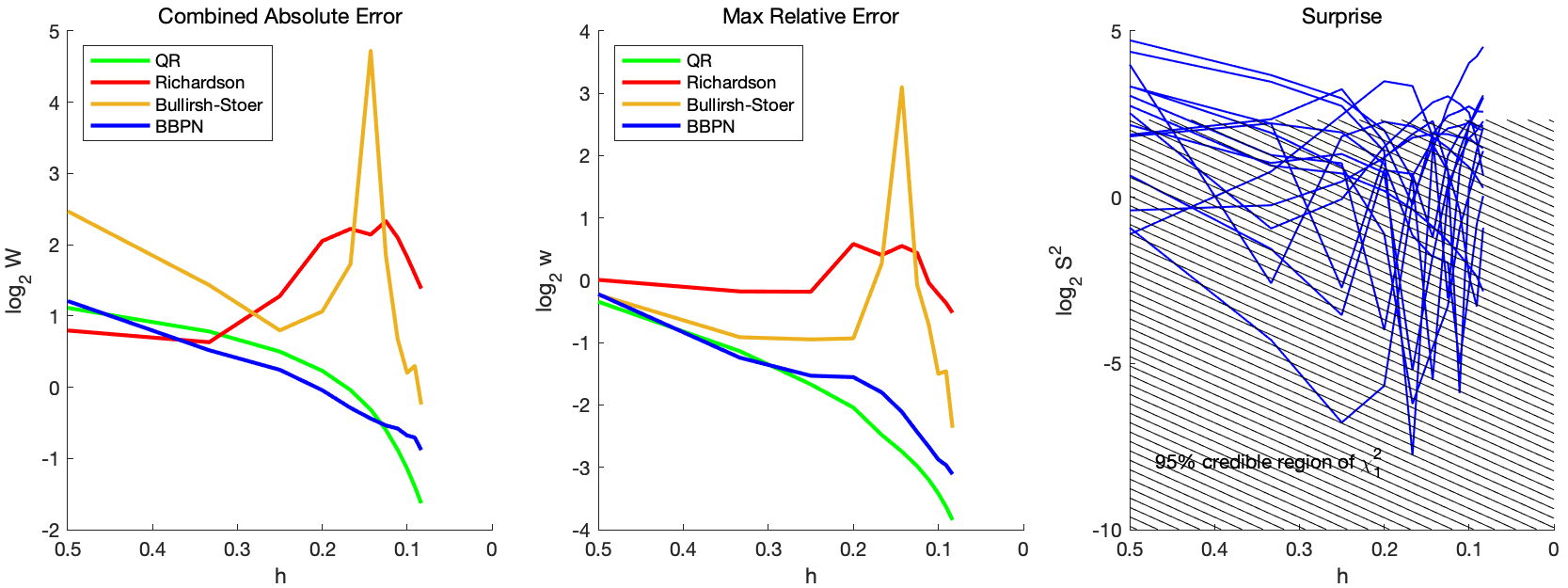}
    \caption{Eigenvalue Problems:
    Left: the combined absolute error $W$ for the classic QR algorithm (green), the traditional extrapolation methods of Richardson (red) and Bulirsh--Stoer (yellow) using the data obtained in the run of the QR algorithm, and \ac{BBPN} (blue). Centre: the maximum relative error $w$ for the same methods. 
   Right: the surprise $S$ of \emph{individual} eigenvalues of the $15 \times 15$ matrix, plotted over the 95\% central probability region of a $\chi_1^2$ random variable.
    }
    \label{fig:eig comparison}
\end{figure}

In a similar manner to \Cref{fig:ode figure}, we plot in the left-hand pane of \Cref{fig:eig comparison} the (log-) error $W$ for several methods \textemdash the classic QR algorithm in green, then the traditional extrapolation methods of Richardson and Bulirsch--Stoer (using the data obtained in the run of the QR algorithm) in red and yellow respectively. From the definition of $W$, this `combined absolute error' is formed by considering the norm of the error vector of \emph{all} eigenvalues. (The centre pane gives the (log-) maximum relative error $w$, \emph{i.e.} $\max_i [(\hat{\lambda}_i - \lambda_i) / \lambda_i]$, where $\hat{\lambda}$ is the vector of true eigenvalues, and is provided since this is a more familiar presentation of error in eigenvalue problems in numerical analysis.) 

It is seen that polynomial and even rational function interpolation are not robust in this setting, and give errors significantly larger than simply the most accurate single QR-produced estimate.
BBPN does not suffer the same issue, possibly because the nonparametric interpolant has favourable stability properties, and it is somewhat competitive with the traditional QR algorithm, at the cost of additional computation but with the additional richness of output that a \ac{PN} method provides.

The right-hand pane shows the (log-squared-) surprise of \emph{individual} eigenvalues of the $15 \times 15$ matrix, plotted over the 95\% central probability region of a $\chi_1^2$ random variable. 
This shows that the predictions provided for the majority of the 15 eigenvalues are well-calibrated, but that a small number of predictions are overconfident.
This is a promising early result for a problem with no previous \ac{PN} method in existence, as well as one in which $\alpha$ has to be inferred due to the absence of a canonical parameterisation for $h$; see \Cref{app: stationarity}.

\subsection{Kuramoto–Sivashinsky Equation (KSE)}
\label{app: PDEs}

In this section we provide further detail for the \ac{PDE} problem presented in \Cref{sec:PDEs}.

Numerical solutions to the \ac{KSE} were computed on the spatial grid $x\in\{0,0.001,0.002,\dots, 1\}$ and over time segments $t_{i,j} = jh_{i}$ for $j\in\{0,1,\dots, m_j\}$, where $m_j = \lfloor{200/h_{i}}\rceil$ with $\lfloor\bullet\rceil$ denoting the nearest integer function and $h_j$ the time-step parameter. After transformation into Fourier space, solutions were computed using a fourth-order Runge--Kutta numerical integrator ETDRK4 \cite{Kassam2005}.   

\subsubsection{Fourier Transform to Employ the ETDRK4 Numerical Integration Scheme} We discretise the spatial domain using a Fourier spectral transformation. That is, we set
\[
u(x,t) \approx \sum_{k\in\Omega_{k}} \tilde{u}_{k}(t)\exp^{ikx/L},
\]
in (\ref{eq:KSE}), where $\Omega_{k}$ denotes the set of wave-numbers. Doing so returns the Fourier transformed \ac{KSE},
\begin{equation}\label{eq:fourier}
\frac{d}{dt}\tilde{u}_{k}(t) + \left(\frac{k^4}{L^4}-\frac{k^2}{L^2}\right)\tilde{u}_{k}(t) + \frac{ik}{2L}\tilde{v}_{k}(t) = 0, \qquad t>0,
\end{equation}
where
\[
\tilde{v}_{k}(t) = \frac{1}{2\pi L}\int_{-\pi L}^{\pi L}u^{2}(x,t)\exp^{-ikx/L}dx\approx \frac{1}{N}\sum_{l =0}^{N-1}u^{2}(x_{l},t)\exp^{-ikx_{l}/L}
\]
with $N = 1/\delta{x}$ and $\delta{x}$ denoting the spatial step-size, and on assuming that both the solution and spatial derivative are periodic in $x$, \emph{i.e.}, 
\[
u(x,t) = u(x + 2\pi L, t) \quad\text{and}\quad \frac{\partial}{\partial x}u(x,t) = \frac{\partial}{\partial x}u(x + 2\pi L, t), \qquad t\geq 0,
\]
for some user defined length scale $L$ (which we take to be $L=1/2\pi$ in our simulation). See \cite{Kassam2005} for a complete description of the fourth-order ETDRK4 scheme, as well as example \texttt{MATLAB} code used to compute solutions to \eqref{eq:fourier}.

\end{document}